\documentclass[12pt]{article}

\usepackage{ifpdf}
\ifpdf 
    \usepackage[pdftex]{graphicx}   
    \usepackage[pdftex,     
            plainpages=false,   
            breaklinks=true,    
            colorlinks=true,
            pdftitle=My Document
            pdfauthor=My Good Self
           ]{hyperref} 
    \usepackage{thumbpdf}
\else 
    \usepackage{graphicx}       
    \usepackage{hyperref}       
\fi 
\usepackage[utf8]{inputenc}
\usepackage{amsmath}
\usepackage{amsthm}
\usepackage{amsfonts}
\usepackage[english]{babel}
\usepackage{rotating}
\usepackage{a4wide}
\usepackage{amssymb}
\usepackage{enumerate}
\usepackage{stackrel}
\usepackage{hyperref}
\usepackage{framed, color}
\usepackage{color}
\usepackage{authblk}
\usepackage{dsfont}

\newcommand{\la}{\lambda}
\renewcommand{\L}{\mathcal{L}}
\renewcommand{\b}{\beta}
\renewcommand{\a}{\alpha}
\newcommand{\n}{\nabla}
\newcommand{\R}{\mathbb{R}}
\newcommand{\p}{\partial}
\newcommand{\ep}{\epsilon}
\renewcommand{\t}{\theta}
\renewcommand{\o}{\omega}
\newcommand{\ds}{\displaystyle}

\newcommand{\KS}{\mathcal{KS}_n}

\def\be{\begin{equation}}
\def\ee{\end{equation}}
\def\bea{\begin{eqnarray}}
\def\eea{\end{eqnarray}}
\def\bean{\begin{eqnarray*}}
\def\eean{\end{eqnarray*}}

\newtheorem{theorem}{Theorem}
\newtheorem{defn}{Definition}
\newtheorem{propo}{Proposition}
\newtheorem{oq}{Open Question}
\newtheorem{lema}{Lemma}
\newtheorem{coro}{Corollary}
\newtheorem{conj}{Conjecture}

\newcommand{\uwidehat}[1]{%
  \mathpalette\douwidehat{#1}%
}

\makeatletter
\newcommand{\douwidehat}[2]{%
  \sbox0{$\m@th#1\widehat{\hphantom{#2}}$}%
  \sbox2{$\m@th#1x$}
  \sbox4{$\m@th#1#2$}
  \dimen0=\ht0
  \advance\dimen0 -.8\ht2
  \dimen2=\dp4
  \rlap{%
    \raisebox{\dimexpr\dimen0-\dimen2}{%
      \scalebox{1}[-1]{\box0}%
    }%
  }%
  {#2}%
}
\makeatother

\title{Semi-Riemannian manifolds with linear differential conditions on the curvature}
\author{José M.M. Senovilla$^{1,2}$}
\affil{$^1$ Departamento de F\'{\i}sica, Facultad de Ciencia y Tecnolog\'ia, Universidad del Pa\'{\i}s Vasco UPV/EHU, Apartado 644, 48080 Bilbao, Spain.\\
$^2$ EHU Quantum Center, Universidad del Pa\'{\i}s Vasco UPV/EHU, Bilbao, Spain.}
\date{}

\begin{document}
\maketitle

\begin{abstract}
Semi-Riemannian manifolds that satisfy (homogeneous) linear differential conditions of arbitrary order $r$ on the curvature are analyzed. They include, in particular, the spaces with ($r^{th}$-order) recurrent curvature, ($r^{th}$-order) symmetric spaces, as well as entire new families of semi-Riemannian manifolds rarely, or never, considered before in the literature --such as the spaces whose derivative of the Riemann tensor field is recurrent, among many others. Definite proof that all types of such spaces do exist is provided by exhibiting explicit examples of all possibilities in all signatures, {\em except in the Riemannian case with a positive definite metric}. Several techniques of independent interest are collected and presented. Of special relevance is the case of Lorentzian manifolds, due to its connection to the physics of the gravitational field. This connection is discussed with particular emphasis on Gauss-Bonnet gravity and in relation with Penrose limits. Many new lines of research open up and a handful of conjectures, based on the results found hitherto, is put forward.

\vspace{1cm} 

\end{abstract}
\tableofcontents

\section{Introduction}\label{sec1}
The investigation of a type of tensor fields in Lorentzian manifolds named ``super-energy tensors'' \cite{S}  led to the following very simple question: if $R$ is the curvature tensor, which Lorentzian manifolds satisfy 
$$
\n^2 R=0 \, ?
$$
 Traditional results in Riemannian geometry \cite{TA,NO} stated that $\n^2 R=0$ implied local symmetry: $\n R=0$, and this seems to be the reason why a higher-order condition such as $\n^2 R=0$, or others that one can imagine, had not been explored before in  non-Riemannian possibilities. However, such a ``negative'' result does not hold in the general semi-Riemannian case. 
 
In the Lorentzian case, the main ideas were presented in \cite{SN1}, where the fundamental result of the existence of a parallel null vector field was proven, which led to a full resolution of the problem first in four dimensions \cite{BSS1}, and then in the general case in \cite{BSS2,BSS3,Tesis} and , independently, in \cite{AG}.

These achievements lead, somehow naturally, to the study of semi-Riemannian manifolds that satisfy a homogeneous linear differential condition on the curvature:
\begin{equation}\label{LDC}
 \n^{r}R+t^{(1)}\otimes \n^{r-1}R+t^{(2)}\otimes \n^{r-2}R+\ldots+t^{(r-1)}\otimes \n R+t^{(r)}\otimes R=0  
 \end{equation} 
 for some $m$-covariant tensor fields $t^{(m)}$, $m\in\{1,\dots,r\}$.
 
This family of manifolds was presented in two invited talks \cite{Jaen,SN2} at meetings on Lorentzian geometry and General Relativity. After this, a very important paper from 1985 by Kaigorodov \cite{Kai} came to my knowledge
where the semi-Riemannian manifolds satisfying \eqref{LDC} where studied in some detail --but not so those with $\n^2 R=0$. 
Curiously enough, the results from \cite{Kai} and \cite{Jaen,SN2} were mostly complementary or independent, and the techniques used different. The purpose of this paper is to publish the main results, the basic ideas and some important open questions combining partly the results in \cite{Kai} with those in \cite{Jaen,SN2} which are published herein for the first time. I firmly believe that this line of research will be helpful and interesting in the areas of semi-Riemannian geometry and in the physics of the gravitational field.

 This paper is structured as follows: in Section \ref{sec2} the family of semi-Riemannian manifolds satisfying \eqref{LDC} is introduced, their main properties listed and a classification presented and analyzed. Section \ref{exp} presents a model family of explicit semi-Riemannian manifolds (Theorem \ref{exa}) satisfying \eqref{LDC} for all possible $r$ and for all signatures {\em except} the positive-definite one --because they possess a parallel lightlike vector field. In Section \ref{sec:cases} all known results for particular values of $r$ and of the tensor fields $t^{(m)}$ in \eqref{LDC} are collected. Section \ref{sec5} is devoted to a brief analysis of the physical relevance of the Lorentzian manifolds satisfying \eqref{LDC} , in particular their connection with the Penrose limits \cite{Pen} and its possible relevance in higher dimensional theories such as Gauss-Bonnet gravity \cite{FCCM} or those derived from string theory. Finally, section \ref{sec6} contains a discussion with the open lines of research and presents a set of conjectures that seem to be valid taking all the known results into account.
 
This is complemented with three appendices, with contents of independent interest, collecting some interesting and classical techniques that are helpful in proving the main results, such as the existence of ``generic points'', the use of particular homothetic vector fields that are gradients and arise naturally in this context, and the general results concerning when a tensor field $T$ with vanishing $2^{nd}$-order covariant derivative is actually parallel.

Before that, the next subsection fixes the notation and the terminology used throughout. This is basically standard, and thus it can be skipped in a first reading, or if the reader is familiar with the subject. Analogously, it may be advisable to only look at the Appendices when needed upon reading the main text. 

\subsection{Notation and terminology}\label{subsec:notation}
Let $(M,g)$ be an $n$-dimensional semi-Riemannian manifold of arbitrary signature, $\nabla$ its Levi-Civita connection and $R(X,Y)Z=\n_X\n_YZ-\n_Y\n_XZ-\n_{[X,Y]}Z$ its curvature tensor field. The signature is Riemannian if the metric is positive definite, and is called Lorentzian if the signature is $(1,n-1)$ --or $(n-1,1)$. The use of abstract index notation will be extremely convenient in many cases in this work, but it will not refer to the use of any particular basis or coordinate choice --under explicitly stated on the contrary. Thus, the Riemann tensor will be denoted by $R$ or, in abstract index form, by
$R^{\a}\,_{\b\la\mu}$ such that
$$(\nabla_\la \nabla_\mu -\nabla_\mu\nabla_\la) X^\a = R^\a{}_{\b\la\mu} X^\b.$$
This is sometimes called the Ricci identity. $Ric$ will denote the corresponding Ricci tensor or, in abstract index notation $R_{\a\b}:=R^{\rho}\,_{\a\rho\b}$ while $S$ denotes the scalar curvature. When using indices the metric is $g_{\a\b}$ and its contravariant version (or its inverse) is $g^{\a\b}$. In index notation any vector field $X$ carries a contravariant index  $X^\a$, and its metrically equivalent one-form will be denoted by $\uwidehat{X}$ or, in abstract indices, by $X_\a:= g_{\a\b}X^\b$. Thus, the covariant Riemann tensor is $R_{\a\b\la\mu}=g_{\a\rho}R^\rho\,_{\b\la\mu}$. Here and throughout the Einstein summation convention will be used. When using index notation, round brackets will mean symmetrization in the enclosed indices, and square brackets anti-symmetrization.

Within the class of metrics that satisfy \eqref{LDC} there are obvious classic examples for some particular values of the tensor fields $t^{(m)}$, $m=1,\dots ,r$ which already have their nomenclature. Outstanding examples are the first order conditions $ (r=1) $ corresponding to {\em recurrent} spaces if $t^{(1)}\neq 0$, and to {\em locally symmetric} manifolds if $t^{(1)}= 0$. More generally, $r^{th}$-order recurrent --or simply $r$-recurrent-- spaces are those with $t^{(1)}=\dots =t^{(r-1)}=0$ and $t^{(r)}\neq 0$, while $r^{th}$-order symmetric --or simple $r$-symmetric-- spaces are defined by $t^{(m)}=0$ for all $m\in\{1,\dots ,r\}$.

Recall that a semi-Riemannian manifold is called reducible when the holonomy group leaves a non-trivial subspace of each tangent space invariant; further, it is non-degenerately reducible if it leaves a non-degenerate subspace—that is, such that the restriction of the metric is non-degenerate—invariant. In the latter case the manifold is (locally)  {\em decomposable} in the sense that $(M,g)$ is a direct product of two semi-Riemannian manifolds $(M_1,g_1)$ and $(M_2,g_2)$ with $M=M_1\times M_2$ and $g=g_1\oplus g_2$.
The following terminology will be used:
\begin{defn}[$r^{th}$-symmetric extensions]\label{ext}
An {\rm $r$-symmetric extension of a semi-Riemannian manifold $(M,g)$ } is a direct product manifold that is decomposable into $(M,g)$ and a $r^{th}$-order symmetric semi-Riemannian manifold $(\overline{M},g^{\overline{M}})$. 

If $r=0$,  $0$-symmetric extensions will be simply called {\rm flat extensions} of the manifold, meaning that $(\overline{M},g^{\overline{M}})$ is a flat semi-Riemannian manifold, and if $r=1$,  $1^{st}$-order symmetric extensions will be referred to as {\rm locally symmetric extensions}.\end{defn}

It is important to point out that there exist semi-Riemannian manifolds characterized by conditions {\em non-linear} on the curvature. The most relevant such condition is that of semi-symmetry \cite{Cartan,Szabo,Szabo1}, which was actually introduced due to the absence of $r$-symmetric and $r$-recurrent spaces in Riemannian geometry, see section \ref{sec:cases}. This is generalized to higher non-linearities as follows \cite{Kai}
\begin{defn}[$\frac{1}{2p}$-symmetric semi-Riemannian manifolds]
   A semi-Riemannian manifold $(M,g)$ is called {\em $\frac{1}{2p}$-symmetric} with $p\in \mathbb{N}$  if it satisfies the condition  $$\n_{[\la_1}\n_{\la_2]}\n_{[\la_3}\n_{\la_4]}\ldots \n_{[\la_{2p-1}}\n_{\la_{2p}]} R^{\a}\,_{\b\la\mu}=0.$$
   In particular,  if $p=1$ it is called a {\em semi-symmetric space} and the condition becomes $R(X,Y) \cdot R=0$, which is equivalent to $\n_{[\la}\n_{\mu]} R^{\a}\,_{\b\gamma\delta}=0$.
\end{defn} 
Here $R(X,Y) \cdot R=0$ is a shorthand \cite{Szabo} for the following expression in abstract indices
\be
R^{\rho}{}_{\alpha\lambda\mu}R_{\rho\beta\gamma\delta}+
R^{\rho}{}_{\beta\lambda\mu}R_{\alpha\rho\gamma\delta}+
R^{\rho}{}_{\gamma\lambda\mu}R_{\alpha\beta\rho\delta}+
R^{\rho}{}_{\delta\lambda\mu}R_{\alpha\beta\gamma\rho}=0 .\label{RicR}
\ee

\section{$\KS$ spaces: definition and properties}\label{sec2}
The main subject of this paper is the following set of semi-Riemannian manifolds:
  \begin{defn}[$\KS$ spaces. \cite{KaiA,Kai,Jaen,SN2}]\label{defKS}
    An $n$-dimensional semi-Riemannian manifold $(M,g)$ {\em satisfies a linear differential condition of order $r\in \mathbb{N}$ on the curvature tensor $R$ associated to a \underline{fixed} tensor field $B\in T^1_{r+3}(M)$}  if there exist $r$ covariant tensor fields $t^{(m)}\in T^0_{m}(M)$ with $m=1,\ldots,r$, such that 
  \begin{equation}\label{nHLDC}
   \n^{r}R+t^{(1)}\otimes \n^{r-1}R+t^{(2)}\otimes \n^{r-2}R+\ldots+t^{(r-1)}\otimes \n R+t^{(r)}\otimes R=B.  
    \end{equation}
 If $B=0$, then $(M,g)$, called $\KS$, satisfies the {\em homogeneous} linear differential condition of order $r$ on the curvature \eqref{LDC}.
      \end{defn}
Some research has been done in the non-homogeneous case for particular manifolds (see for example \cite{nh1}-\cite{nh4}). In section \ref{exp} some explicit examples with $B\neq 0$ will also be given. Nevertheless, this work will be dealt mainly with the homogeneous case, Eq.\eqref{LDC}, which was first studied by Kaigorodov \cite{KaiA, KaiB, Kai}  and much later re-discovered in \cite{Jaen,SN2}.
 
Despite the linearity of the condition \eqref{LDC} on $R$, there is no superposition principle in the solutions of these equations in the sense that, fixed $M$, the linear combinations of curvature tensors satisfying \eqref{LDC} associated to different metrics on $M$ will not be another solution in general, since both $R$ and $\n$ depend also on the chosen metric. 

Making some abuse of notation and removing the superscript of the tensor fields $t^{(m)}$ in  equation \eqref{LDC} --since they can be identified by their number of subindices in abstract index form--, the equation can be written in the abstract index form as:
 \begin{align}\label{LDCai}
   \n_{\la_1}\ldots\n_{\la_r}R^\a\,_{\b\mu\nu}+t_{\la_1} \n_{\la_2}\ldots\n_{\la_r}R^\a\,_{\b\mu\nu}+\ldots+t_{\la_1\ldots\la_{r-1}} \n_{\la_r} R^\a\,_{\b\mu\nu}+t_{\la_1\ldots\la_r} R^\a\,_{\b\mu\nu}&=0.
 \end{align}
Observe that, written in this form, this definition could be generalized by permuting the subscripts to obtain different linear differential conditions on the curvature, which would lead to a much larger class of manifolds. 
 
 It is obvious that if \eqref{LDC} holds, then linear differential relations with exactly the same structure, and the same set of tensors $t^{(m)}$, are satisfied by the Ricci tensor and the scalar curvature:
\bea 
 \n^{r}Ric +t^{(1)}\otimes \n^{r-1}Ric+t^{(2)}\otimes \n^{r-2}Ric+\ldots+t^{(r-1)}\otimes \n Ric+t^{(r)}\otimes Ric=0,\label{LDCRic}\\
  \n^{r}S+t^{(1)}\otimes \n^{r-1}S+t^{(2)}\otimes \n^{r-2}S+\ldots+t^{(r-1)}\otimes \n S+t^{(r)}\otimes S=0\label{LDCS}. 
\eea 
or with indices
\bean
 \n_{\la_1}\ldots\n_{\la_r}R_{\b\nu}+t_{\la_1} \n_{\la_2}\ldots\n_{\la_r}R_{\b\nu}+\ldots+t_{\la_1\ldots\la_{r-1}} \n_{\la_r} R_{\b\nu}+t_{\la_1\ldots\la_r} R_{\b\nu}=0,\\
  \n_{\la_1}\ldots\n_{\la_r}S+t_{\la_1} \n_{\la_2}\ldots\n_{\la_r}S+\ldots+t_{\la_1\ldots\la_{r-1}} \n_{\la_r} S +t_{\la_1\ldots\la_r} S=0.
\eean
Now, let $C$ denote the Weyl tensor of $(M,g)$. Using \eqref{LDC}, \eqref{LDCRic} and \eqref{LDCS} it is easily verified that the Weyl tensor satisfies the same relation too:
$$
  \n^{r}C+t^{(1)}\otimes \n^{r-1}C+t^{(2)}\otimes \n^{r-2}C+\ldots+t^{(r-1)}\otimes \n C+t^{(r)}\otimes C=0.\label{LDCW}
$$
This opens the  door to define semi-Riemannian manifolds that comply just with  this relation for $C$, or with \eqref{LDCRic}, but not for the full $R$.

If \eqref{LDC} holds, by computing the $r^{th}$-order derivative along any geodesic on $(M,g)$ of $R$ (or $C$, $Ric$ or $S$) and using \eqref{LDC} (or \eqref{LDCW}, \eqref{LDCRic} or \eqref{LDCS}) the following result is derived:
\begin{lema}\label{ZeroD}
On a $\KS$ semi-Riemannian manifold $(M,g)$:
\begin{enumerate}
\item if $R$ and its covariant derivatives up to order $(r-1)$ vanish at a given point $x\in M$, then the curvature vanishes everywhere, that is, $(M,g)$ is flat.
\item if the Weyl tensor $C$ and its covariant derivatives up to order $(r-1)$ vanish at a given point $x\in M$, then $C$ vanishes everywhere and $(M,g)$ is conformally flat.
\item if $Ric$ and its covariant derivatives up to order $(r-1)$ vanish at a given point $x\in M$, then the Ricci tensor vanishes everywhere, that is, $(M,g)$ is Ricci-flat.
\item if $S$ and its covariant derivatives up to order $(r-1)$ vanish at a given point $x\in M$, then the scalar curvature vanishes everywhere.
\end{enumerate}
\end{lema}
  
Notice that if $(M,g)$ satisfies equation \eqref{LDC} for some $r$ then, by taking successive covariant derivatives there, it will also satisfy the same type of equation for all $s>r$ with some properly defined new covariant tensor fields $\overline{t}\,^{(m)}$ with $m=1,\ldots, s$. 
Hence,  the following definition is necessary:
  \begin{defn}
     A $\KS$ semi-Riemannian manifold is said to be {\em proper of order $r$} if $r$ is the minimum value for which the manifold satisfies a condition of  type \eqref{LDC}.
       \end{defn}

\subsection{Non-decomposability}\label{sec31}
In the theory of semi-Riemannian manifolds $(M,g)$ the question of decomposability, or non-degenerate reducibility, is of paramount importance
(see \cite{NO,RH} for the Riemannian case and \cite{WU} for the semi-Riemannian case). For $\KS$ semi-Riemannian manifolds, and
using Definition \ref{ext}, the following powerful theorem, which is a broad generalization of Theorem 2.1 in \cite{WAL}, can be proven:
\begin{theorem}\label{red}\cite{SN2,Jaen}
  If $(M,g)$ is a non-flat proper order $r$ $\KS$ semi-Riemannian manifold, then either it is  $r^{th}$-symmetric, or it is not decomposable except for $k^{th}$-symmetric extensions with $k\leq r$.
    
Furthermore, if $\overline{m}$ is the greatest interger in  $\{1,\ldots,r\}$ such that $t^{(\overline{m})}\neq 0$, then only $s^{th}$-symmetric extensions are allowed, for $s\leq r-\overline{m}$. 
  \end{theorem}
{\bf Proof.}
Suppose that $(M,g)$ is reducible, so that there exist two semi-Riemannian manifolds $(M_1,g_1),(M_2,g_2)$ such that $M=M_1\times M_2$ and $g=g_1\oplus  g_2$. Let us prove that then $(M,g)$ is a $k^{th}$-symmetric extension of $(M_1,g_1)$ or of $(M_2,g_2)$. 
Let $A,B,C,\ldots$ be the indices for $(M_1,g_1)$ and $a,b,c,\ldots$ the indices for $(M_2,g_2)$, then not only the metric but also the covariant derivative $\n$ and the curvature tensor $R$ decompose:
$$g=g_{AB}dx^Adx^B+g_{ab}dx^a dx^b,$$
$$\n=\n^{(1)}+\n^{(2)},$$
$$R=R^{(1)}+R^{(2)},$$
where $g_1=g_{AB}dx^Adx^B$, $g_2=g_{ab}dx^a dx^b$, and $\n^{(i)}$ are the Levi Civita connections and $R^{(i)}$ the curvature tensors of $(M_i,g_i)$, for $i=1,2$. Taking indices $\la_r=A$ and $\a\b\mu\nu=abcd$ in the equation \eqref{LDCai}, the $r^{th}$-order homogeneous condition  becomes
$$t_{\la_1\ldots\la_{r-1}A}R^a\,_{bcd}=0.$$
Analogously, taking indices $\la_r=a$ and $\a\b\mu\nu=ABCD$, \eqref{LDCai} becomes
$$t_{\la_1\ldots\la_{r-1}a}R^A\,_{BCD}=0.$$
Since the manifold is not flat by hypothesis, if $t^{(r)}\neq 0$ then either $R^a\,_{bcd}=0$  or  $R^A\,_{BCD}=0$, and the manifold is a flat extension of either $(M_1,g_1)$ or $(M_2,g_2)$.
If otherwise $t^{(r)}=0$, \eqref{LDCai} becomes
   \begin{equation}\label{f3}
     \n_{\la_1}\ldots\n_{\la_r}R^\a\,_{\b\mu\nu}+t_{\la_1} \n_{\la_2}\ldots\n_{\la_r}R^\a\,_{\b\mu\nu}+\ldots+t_{\la_1\ldots\la_{r-1}} \n_{\la_r} R^\a\,_{\b\mu\nu}=0.   \end{equation} 
    Obviously, if $(M,g)$ is locally symmetric, equation \eqref{f3} is trivially satisfied (and then $r$ would be 1). Suppose that the manifold is not locally symmetric. Then, using the same idea as before take now $\la_{r-1}=A$, $\la_r=e$ and $\a\b\mu\nu=abcd$ on the one hand, and $\la_{r-1}=a$, $\la_r=E$ and $\a\b\mu\nu=ABCD$ on the other hand, to obtain from \eqref{f3}: 
$$
t_{\la_1\ldots\la_{r-2}A}\n_e R^a\,_{bcd}=0,\hspace{1cm}
 t_{\la_1\ldots\la_{r-2}a}\n_E R^A\,_{BCD}=0.
 $$
 Then, if  $t^{(r-1)}\neq 0$ and since $\n R\neq 0$ it follows that either $\n_e R^a\,_{bcd}=0$ or $\n_E R^A\,_{BCD}=0$ and the manifold is a locally symmetric extension of either $(M_1,g_1)$ or $(M_2,g_2)$. 
 If otherwise $t^{(r-1)}= 0$,  
 the last term on left-hand side of \eqref{f3} vanishes:
   $$ \n_{\la_1}\ldots\n_{\la_r}R^\a\,_{\b\mu\nu}+t_{\la_1} \n_{\la_2}\ldots\n_{\la_r}R^\a\,_{\b\mu\nu}+\ldots+t_{\la_1\ldots \la_{r-2}}\n_{\la_{r-1}}\n_{\la_r}R^\a\,_{\b\mu\nu}=0.$$
 Continuing with this reasoning we conclude that either the manifold  $(M,g)$ is a proper $r^{th}$-symmetric space, or it admits a $s^{th}$-symmetric extension with $s\leq r-\overline{m}$.      
\begin{flushright}
{$\blacksquare$}\end{flushright}

Therefore, for $\KS$ spaces proper of order $r$ the cases of interest are the indecomposable ones, any others being mere $s^{th}$-symmetric extensions of them: if $t^{(r)}\neq 0$, the manifold is not decomposable except for flat extensions; if $t^{(r)}=0$ but $t^{(r-1)}\neq 0$, it is not decomposable except for locally symmetric or flat extensions; if $t^{(r)}= 0$ and $t^{(r-1)}= 0$ but $t^{(r-2)}\neq 0$, it is not decomposable except for $2^{nd}$-symmetric, locally symmetric or flat extensions, and so on. 

\subsection{Properties of $t^{(r)}$}
For $\KS$ semi-Riemannian manifolds the following property follows at generic points (these points are defined and discussed in Appendix \ref{sec32}):
\begin{lema}\cite{Jaen,SN2}\label{tr=0}
    If $(M,g)$ is a non-flat $\KS$ semi-Riemannian manifold proper of order $r$, then in a neighbourhood of any generic point $t^{(r)}=0$.
 \end{lema}
 
 {\bf Proof.}
 Applying the second Bianchi identity $\n_{[\la_r}R_{\a\b]\mu\nu}=0$ to equation \eqref{LDCai}, one gets that $$t_{\la_1\ldots\la_{r-1}[\la_r}R_{\a\b]\mu\nu}=0.$$
  Let $x\in M$ be generic and $U(x)\subset M$ an open neighbourhood of $x$ with all of its points generic. Then, contracting the previous expression with $R^{-1}\,^{\mu\nu\sigma\rho}$ on $U(x)$ one obtains that $t_{\la_1\ldots\la_{r-1}[\la_r}\delta_{\a]}^\sigma \delta_{\beta}^{\rho}=0$. Contracting $\sigma$ with $\alpha$, and $\rho$ with $\beta$, one sees that
 $$t_{\la_1\ldots\la_r}=0\, \,  \text{on} \, \, \, U(x) .$$
  \begin{flushright}
{$\blacksquare$}\end{flushright}

 \begin{lema}\cite{Kai}\label{p1}
   Let $(M,g)$ be a non-flat semi-Riemannian manifold of type $\KS$ with $r\geq 2$. If $t^{(r-1)}=0$, then $t^{(r)}$ is symmetric in its last two indices. 
 \end{lema}
 {\bf Proof.}
 If $t^{(r)}=0$ the claim is obvius. Thus consider the case with $t^{(r)}\neq 0$ --which automatically implies that we are dealing with non-generic points due to Lemma \ref{tr=0}.
The following identity for the Riemann tensor holds on any semi-Riemannian manifold \cite{WAL}:
$$\n_{[\a_1} \n_{\a_2]}R_{\b_1\b_2\la_1\la_2}+\n_{[\b_1} \n_{\b_2]}R_{\la_1\la_2\a_1\a_2}+\n_{[\la_1} \n_{\la_2]}R_{\a_1\a_2\b_1\b_2}=0.$$
Taking the covariant derivative of this identity $(r-2)$-times and using equation \eqref{LDCai} and the fact that $t^{(r-1)}=0$, one gets
\begin{equation}\label{12p}
  t_{\la_1\ldots\la_{r-2}[\a_1\a_2]} R_{\b_1\b_2\la_1\la_2}+t_{\la_1\ldots\la_{r-2}[\b_1\b_2]}R_{\la_1\la_2\a_1\a_2}+t_{\la_1\ldots\la_{r-2}[\la_1\la_2]}R_{\a_1\a_2\b_1\b_2}=0.
\end{equation}
  Now, define the 2-form $T^{(\Omega)}_{\a\b}=t_{\la_1\ldots\la_{r-2}[\a\b]}\Omega^{\la_1\ldots\la_{r-2}}$, for an arbitrary $(r-2)$-contravariant tensor field $\Omega$, so that from equation \eqref{12p} one derives
  \begin{equation}\label{12po}
 T^{(\Omega)}_{\a_1\a_2} R_{\b_1\b_2\la_1\la_2}+T^{(\Omega)}_{\b_1\b_2}R_{\la_1\la_2\a_1\a_2}+T^{(\Omega)}_{\la_1\la_2}R_{\a_1\a_2\b_1\b_2}=0.
\end{equation} 
By using collective indices for each skew-symmetric pair here, say ${\cal A}=\a_1\a_2$, ${\cal B}=\b_1\b_2$ and ${\cal C}=\la_1 \la_2$, this can be viewed as the vanishing of the completely symmetrized product $T^{(\Omega)}_{({\cal A}}R_{{\cal B}{\cal C})}=0$, which readily implies that either $T^{(\Omega)}=0$ or $R=0$. As $(M,g)$ is not flat by hypothesis, it follows that $T^{(\Omega)}=0$ for all possible $\Omega$, that is, $t_{\la_1\ldots\la_{r-2}[\a\b]}=0$.

%
 \begin{flushright}
  $\blacksquare$
\end{flushright}

\subsection{An interesting property}
Here we prove that the existence of a vector field $X$ such that $\n X =c \mathds{1}$, $c\in \mathbb{R}$, is not possible in general $\KS$ spaces unless $c=0$. The main properties of such vector fields are collected in Appendix \ref{sec33}.
\begin{propo}\cite{Jaen,SN2} \label{hv4}
  Let $(M,g)$ be a semi-Riemannian manifold of type $\KS$ and assume that the manifold admits a vector field $X$ such that $\n X=c \mathds{1}$ with $c\in\R-\{ 0\}$. Then, $(M,g)$ is flat.
\end{propo}
{\bf Proof.}
Contracting \eqref{LDCai} with $X^{\la_2}$ and applying Lemma \ref{hv1} (d) and Lemma \ref{hv3} in Appendix \ref{sec33} one gets:
\begin{align*}-c(r+1)&\n_{\la_1}\n_{\la_3}\ldots\n_{\la_r}R^{\a}\,_{\b\la\mu}+\left(t_{\la_1\rho}X^{\rho}-crt_{\la_1}\right)\n_{\la_3}\ldots\n_{\la_r}R^{\a}\,_{\b\la\mu}\\&+t_{\la_1\rho\la_3}X^\rho\n_{\la_4}\ldots\n_{\la_r}R^{\a}\,_{\b\la\mu}+\ldots+t_{\la_1\rho\ldots\la_r}X^\rho R^{\a}\,_{\b\la\mu}=0,\end{align*}
which is nothing else than another homogeneous linear differential equation of one order less. Contracting now with $X^{\la_3}$ and using the same idea, one gets:
\begin{align*}(-c)^2&(r+1)r\n_{\la_1}\n_{\la_4}\ldots\n_{\la_r}R^{\a}\,_{\b\la\mu}\\&+\left((-c)^2r(r-1)t_{\la_1}-c(r-1)t_{\la_1\rho}X^{\rho}+t_{\la_1\rho_1\rho_2}X^{\rho_1}X^{\rho_2}\right)\n_{\la_4}\ldots\n_{\la_r}R^{\a}\,_{\b\la\mu}\\&+t_{\la_1\rho_1\rho_2\la_4}X^{\rho_1} X^{\rho_2} \n_{\la_5}\ldots\n_{\la_r}R^{\a}\,_{\b\la\mu}+\ldots+t_{\la_1\rho_1\rho_2\la_4\ldots\la_r}X^{\rho_1} X^{\rho_2} R^{\a}\,_{\b\la\mu}=0 .\end{align*}
  Following the same reasoning and contracting with $X^{\la_4}$,\ldots,$X^{\la_r}$ it is possible to conclude that $\n_{\la_1} R^{\a}\,_{\b\la\mu}= \overline{t}_{\la_1}R^{\a}\,_{\b\la\mu}$ for some 1-form $\overline{t}$.
 However, contracting here with $X_\a$ and using \eqref{XR} one obtains
 $$
 X_\a \n_{\la_1} R^{\a}\,_{\b\la\mu}=0 =-R^{\a}\,_{\b\la\mu}\n_{\la_1} X_\a =-c R_{\la_1\b\la\mu}
 $$
 so that the metric is flat.
  \begin{flushright}
{$\blacksquare$}\end{flushright}

\subsection{Classification}
Observe that lowering the index $\a$, skew-symmetrizing the indices $\la_r,\a,\b$ in equation $\eqref{LDCai}$ 
 and using the second Bianchi identity $\n_{[\a}R_{\b\la]\mu\nu}=0$, one gets that 
 \begin{equation}\label{gt}
   t_{\la_1\la_2\ldots[\la_r}R_{\a\b]\la\mu}=0.
 \end{equation}
This statement is vacuous at generic points as a consequence of Lemma \ref{tr=0}. However, whenever $t^{(r)}\neq 0$ this leads to the classification of $\KS$ spaces. We start with 
\begin{propo}\cite{Kai}\label{Lexists}
  Let $(M,g)$ be a $\KS$ space with $t^{(r)}\neq 0$. Then, there exist a non-vanishing one-form $\omega \in \Lambda(M)$ such that 
  \begin{equation}\label{pb}
  \omega_{[\gamma}R_{\a\b]\la\mu}=0. 
  \end{equation}
\end{propo}
{\bf Proof.} This is basically the content of equation \eqref{gt}. For each $(r-1)$-contravariant tensor field $\Omega$ it is possible to define the following one-form associated to $\Omega$:
$$t^{(\Omega)}_\a={\Omega}^{\la_1\ldots\la_{r-1}}t_{\la_1\la_2\ldots\la_{r-1}\a}.$$
At least one of these $t^{(\Omega)}$ does not vanish, as otherwise $t^{(r)}=0$ contrary to the assumption. Choose then an $\Omega$ such that $\omega :=t^{(\Omega)}\neq 0$. Equation \eqref{gt} provides the result.
\begin{flushright}
  $\blacksquare$
\end{flushright}

\begin{defn}[The sub-module $\L$]\label{Lspace}
Let $\L$ denote the set of all $\omega\in \Lambda(M)$ satisfying \eqref{pb}.
\end{defn}
Obviously, by linearity, $\L$ is a sub-module of $\Lambda(M)$. As we see below, the dimension of $\L$ cannot be greater than two, which together with equation \eqref{pb} will lead to a decomposition of the curvature tensor depending on the dimension of $\L$. This result is very powerful since it gives the possibility of completing a classification of $\KS$-spaces. To prove this theorem some previous results are needed. The following lemma is well known.
 \begin{lema}\label{auxlem}
   If there exists a one-form $\omega$ such that  $\omega_{[\gamma}R_{\a\b]\la\mu}=0$, then there exist a 2-covariant symmetric tensor field $Q$ such that $R_{\a\b\la\mu}=4 \omega_{[\a}Q_{\b][\la}\omega_{\mu]}$.
 \end{lema}
The tensor field $Q_{\b\la}$ is not uniquely defined, as there is the gauge freedom
\be\label{gauge}
Q_{\b\la} \longrightarrow Q_{\b\la} +\omega_\b \tau_\la +\omega_\la \tau_\b
\ee
for arbitrary $\tau\in \Lambda(M)$. This freedom can be used to simplify $Q$. Basically, all components of $Q$  along $\omega$ can be removed. Thus, $Q$ --and a fortiori $R$-- contains at most $n(n-1)/2$ independent components. 

\begin{lema}\label{co2}
    If there exists two linearly independent one-forms $\omega, \rho\in \Lambda(M)$ such that $R_{\a\b[\la\mu}\omega_{\nu]}=0$ and $R_{\a\b[\la\mu}\rho_{\nu]}=0$, then $R_{\a\b\la\mu}=4A\, \omega_{[\a}\rho_{\b]}\rho_{[\la}\omega_{\mu]}$ for some scalar $A$.
\end{lema}
 {\bf Proof.}
By Lemma \ref{auxlem} there exists a 2-covariant symmetric tensor field $Q$ such that $R_{\a\b\la\mu}=4 \omega_{[\a}Q_{\b][\la}\omega_{\mu]}$. But since $R_{\a\b[\la\mu}\rho_{\nu]}=0$ also, using the symmetry in the interchange of skew-symmetric pairs for $R$ and the symmetry of $Q$, the only possibility is that $Q_{\b\la}=A\, \rho_\b \rho_\la +\omega_\b \tau_\la +\omega_\la \tau_\b$ for some scalar $A$ and some $\tau\in \Lambda(M)$. But using the gauge freedom \eqref{gauge} we arrive at $Q_{\b\la} =A\, \rho_\b \rho_\la $ as required. 
\begin{flushright}
$\blacksquare$
\end{flushright}
In this situation $R=A F\otimes F$ where $F\in \Lambda^2 M$ is the two-form $F=\omega\wedge \rho$, and $R$ contains a  single independent component.
\begin{lema}[dim$\L \leq 2$]\label{Ldim}
In general, dim$\L$ cannot be greater than 2.
\end{lema}
{\bf Proof.}
 By Lemma \ref{co2} any third one-form $\sigma \in \Lambda(M)$ satisfying $R_{\a\b[\la\mu}\sigma_{\nu]}=0$ will necessarily satisfy $\omega_{[\la}\rho_\mu \sigma_{\nu]}=0$ -- unless $R_{\a\b\la\mu}=0$--, that is $\omega\wedge \rho \wedge \sigma =0$, meaning that $\omega,\rho$ and $\sigma$ are linearly dependent.
 \begin{flushright}
  $\blacksquare$
\end{flushright}
There is a strong similarity between the sub-module $\L$ and the so-called ``Olszak distribution'' \cite{Ol,DR} defined in conformally recurrent manifolds for the Weyl tensor instead of $R$. Actually, their algebraic properties are identical, however, the Olszak distribution is a 2- or 1-dimensional {\em parallel} distribution, while $\L$ will not have the parallel property in general, at least not in principle.

From Proposition \ref{Lexists} and Lemma \ref{Ldim} the next theorem follows at once.
 \begin{theorem}\label{TeoL}
   If $(M,g)$ is a non-flat $\KS$ space with $t^{(r)}\neq 0$, then {\rm dim}\,$\L\in \{1,2\}$. 
 \end{theorem}

Thus, the following classification is in order:
\begin{defn}\label{classification}
$\KS$ spaces can be classified into the following disjoint classes: Type {\rm 0}, Type {\rm I} and Type {\rm II} according to whether {\rm dim}\,$\L =0, 1$ or $2$ respectively.
\end{defn}

All $\KS$ with $t^{(r)}\neq 0$ are of types {\rm I} or {\rm II} due to Proposition \ref{Lexists}, so that ${\cal L} $ is non-empty. The Riemann tensor assumes the form given in Lemma \ref{auxlem} for Type {\rm I} and that of Lemma \ref{co2}  for Type {\rm II}. However, $\KS$ with $t^{(r)}=0$ may also belong to types {\rm I} or {\rm II} if they happen to have dim $\L >0$, and in this situation the previous forms of the Riemann tensor also apply. On the other hand, 
Type {\rm 0}  requires $t^{(r)}=0$ and all possible manifolds with generic points belong here --at least in the domains with such points. At this stage it is not clear whether there are any possible manifolds of Type 0, and it may turn out that they are impossible. In any case, the outstanding possibility of $\KS$ with $t^{(r)}=0$ are the $r$-symmetric semi-Riemannian manifolds defined by $\n^r R=0$. All known results concerning $r$-symmetry are collected in section \ref{sec:cases}, and all the known cases will belong to types {\rm I} or {\rm II}.

 \begin{propo}\cite{Kai}
  Let $(M,g)$ be a $\KS$ space of type {\rm I} or {\rm II}. Then, the following relations are satisfied:
  \begin{eqnarray}
    \omega_\rho R^\rho\,_{\a\b\la}=2R_{\a[\la}\omega_{\b]}, \,\,\,\,\,\,\,\, \forall \omega\in \L\label{e1}\\
    \omega^\rho R_{\rho\a}=\frac{1}{2}S\omega_\a, \,\,\,\,\,\,\,\, \forall \omega\in \L\label{e2},\\
    {\cal G}:=R^{\a\b\la\mu}R_{\a\b\la\mu}-4R^{\a\b}R_{\a\b}+S^2=0.\label{e3}
  \end{eqnarray}
\end{propo}
 {\bf Proof.}
Take equation \eqref{pb} and contract the indices $\mu$ and $\gamma$ to get \eqref{e1}. Next, contract $\a$ and $\la$ in \eqref{e1} to get \eqref{e2}. Finally, contract \eqref{pb} with $R^{\a\b\la\mu}$ and use first \eqref{e1} once, and then \eqref{e2} twice to obtain \eqref{e3}.
 \begin{flushright}
  $\blacksquare$
\end{flushright}
In \eqref{e3} one recognizes the Gauss-Bonnet scalar ${\cal G}$. This is identically zero if $n\leq 3$, and in $n=4$ has special topological properties ---see section \ref{sec5}. 
\begin{coro} \label{coro2}
   If  $(M,g)$ is a $\KS$ space of type {\rm I} or {\rm II}, then
   \begin{enumerate}[(a)]
   \item The Gauss-Bonnet scalar vanishes
     \item the scalar curvature vanishes if and only if $R_{\a\b}\omega^\b=0,\, \forall \omega\in \L$;
     \item $\omega$ is, at each point, an eigenvector of the Ricci tensor with eigenvalue $\frac{S}{2}$.
   \end{enumerate}
\end{coro}

\begin{coro}\label{Ricflat}
  If  $(M,g)$ is a $\KS$ space of type {\rm I} or {\rm II}, and is also an Einstein space of dimension $n>2$, then it is in fact Ricci flat.
\end{coro}
 {\bf Proof.}
Substituting $R_{\a\mu}=\frac{1}{n} S g_{\a\mu}$ into \eqref{e2} one gets $\frac{1}{n} S=\frac{1}{2}S$ which implies $S=0$ if $n> 2$, and the space is Ricci flat.
 \begin{flushright}
  $\blacksquare$
\end{flushright}
 \begin{propo}\label{NoEinstein}
 If  $(M,g)$ is a $\KS$ Einstein space of type {\rm I} or {\rm II}, then 
 $$
 R^{\a\b\la\mu}R_{\a\b\la\mu}=0.
 $$
 In particular, if $(M,g)$ is Riemannian, it is actually flat.
 \end{propo}
 {\bf Proof.}
 From the vanishing of the Gauss-Bonnet term \eqref{e3} and the previous Corollary \ref{Ricflat} one readily gets $R^{\a\b\la\mu}R_{\a\b\la\mu}=0$. Of course, if the metric is positive definite, this implies the vanishing of $R$.
  \begin{flushright}
  $\blacksquare$
\end{flushright}

\subsubsection{Type {\rm I}}
 In this class of manifolds $\text{dim }\L=1$ and there exists a one-form $\omega\in \Lambda(M)$ such that $\L=<\omega>$. By Lemma \ref{auxlem}, there also exists a 2-covariant symmetric tensor field $Q_{\a\b}$ such that the curvature tensor, the Ricci tensor and the scalar curvature have the form 
 \begin{eqnarray}
 R_{\a\b\la\mu}&=&4\omega_{[\a}Q_{\b][\la}\omega_{\mu]};\label{ctb}\\
 R_{\a\b}&=&\omega_{\a} Q_{\b\rho}\omega^\rho+\omega_{\b} Q_{\a\rho}\omega^\rho-\omega_{\rho}\omega^{\rho}Q_{\a\b}-Q_\rho\,^\rho\omega_\a \omega_ \b;\label{rtb}\\
S&=&2Q_{\rho\sigma}\omega^\rho\omega^\sigma -2\omega^\rho \omega_\rho Q_\sigma\,^\sigma .\label{ectb}
 \end{eqnarray}
 Observe that, if $\overline{X}$ is a vector field such that $\omega(\overline X)=1$, then
 $$Q_{\a\b}=R_{\a\rho\b\sigma}\overline{X}^\rho\overline{X}^\sigma;\,\,\,\,\,\,\,\,\,\, Q=R_{\rho\sigma}\overline{X}^\rho\overline{X}^\sigma.$$
 
 Type {\rm I} can be split into two cases depending on whether or not $\L$ is degenerate. Thus we have
  \begin{itemize}
   \item Type {\rm I}$_\epsilon$ if $\L$ is non-lightlike, that is if $\omega$ is non-null.
   \item Type {\rm I}$_N$ if $\L$ is lightlike, i.e., if $\omega$ is a lightlike one-form.
  \end{itemize}
Clearly, there are no $KS_n$ spaces of type {\rm I}$_N$ in the Riemannian case.
 
Let $(M,g)$ be a $\KS$ of type ${\rm I}_\epsilon$ and set $\epsilon:=g(\omega,\omega)=\pm 1$. Then, using the freedom \eqref{gauge} one can always choose $Q_{\a\b}$ totally orthogonal to $\omega$, $Q_{\a\b}\omega^\b=0$, 
as we do from now on for type ${\rm I}_\epsilon$. Under this condition, equations \eqref{rtb}-\eqref{ectb} reduce to
 \begin{eqnarray}
 R_{\a\b}&=&-\ep Q_{\a\b}-Q_\rho\,^\rho \omega_\a \omega_ \b;\label{rtb1}\\
S&=&-2\ep Q_\rho\,^\rho,\label{ectb1}
 \end{eqnarray}
 Therefore, the 2-covariant tensor field $Q_{\a\b}$ can be written in terms of the Ricci tensor and the scalar curvature as
 $$Q_{\a\b}=\frac{S}{2}\omega_\a \omega_\b-\ep R_{\a\b};\,\,\,\,\,\,\,\,\,\, Q_{\a\b}\omega^\b=0.$$
  Using this form of $Q_{\a\b}$ it is straightforward to derive the next
 \begin{coro}
 The curvature tensor of type ${\rm I}_\ep$ can be written in terms of the Ricci tensor as
   $$R_{\a\b\la\mu}=4\omega_{[\a}{R}_{\b][\la}\omega_{\mu]}.$$
 \end{coro}

For type ${\rm I}_\ep$, all possible decompositions of the Ricci tensor and the metric were given in \cite{KaiB} for the Riemannian case, and for the Lorentzian signature in \cite{Kai}.
 
 \begin{theorem}\label{p5}
    A $\KS$ semi-Riemannian manifold $(M,g)$ of type ${\rm I}_\ep$ is not $\frac{1}{2p}$-symmetric for any value of $p\in \mathbb{N}$.
 \end{theorem}

{\bf Proof.}
We will use the general expression \eqref{ctb} for the curvature tensor with $Q_{\a\b}\omega^\b=0$, as we have already argued. Then, using \eqref{ctb} and the Ricci identity it is possible to prove that
\begin{equation}\label{p=1}
  \n_{[\a_2}\n_{\a_1]}R_{\a\b\la\mu}=2 \left(\omega_{[\la}Q_{\mu][\a}Q_{\b][\a_2}\omega_{\a_1]}+\omega_{[\a}Q_{\b][\la}Q_{\mu][\a_2}\omega_{\a_1]}\right).
\end{equation} 
A somewhat longer calculation, using the above expression and the Ricci identity again,  allows one to prove that
\begin{align}
  \n_{[\a_4}\n_{\a_3]}\n_{[\a_2}\n_{\a_1]}&R_{\a\b\la\mu}=4 \left(\omega_{[\a_3}Q_{\a_4][\la}Q_{\mu][\a}Q_{\b][\a_2}\omega_{\a_1]}+\omega_{[\a_3}Q_{\a_4][\a}Q_{\b][\la}Q_{\mu][\a_2}\omega_{\a_1]}\right.\nonumber\\
  &+\left. \omega_{[\a_3}Q_{\a_4][\a_2}Q_{\a_1][\a}Q_{\b][\la}\omega_{\mu]}+\omega_{[\a_3}Q_{\a_4][\a_2}Q_{\a_1][\la}Q_{\mu][\a}\omega_{\b]}\right)\nonumber \\
  &+4\ep\left(\omega_{\mu} Q^{\rho}\,_{[\a_4}\omega_{\a_3]}\omega_{[\a}Q_{\b][\rho}Q_{\la][\a_2}\omega_{\a_1]}-\omega_{\la} Q^{\rho}\,_{[\a_4}\omega_{\a_3]}\omega_{[\a}Q_{\b][\rho}Q_{\mu][\a_2}\omega_{\a_1]}\right.\nonumber\\
  &\left.+\omega_{\b} Q^{\rho}\,_{[\a_4}\omega_{\a_3]}\omega_{[\la}Q_{\mu][\rho}Q_{\a][\a_2}\omega_{\a_1]}-\omega_{\a} Q^{\rho}\,_{[\a_4}\omega_{\a_3]}\omega_{[\la}Q_{\mu][\rho}Q_{\b][\a_2}\omega_{\a_1]}\right)\label{p=2}
 \end{align}
Contracting \eqref{p=1} with $\omega^\mu$ and $\omega^{\a_2}$ and using the condition of semi-symmetry $(p=1)$, and contracting \eqref{p=2} with $\omega^\mu$ and $\omega^{\a_4}$ and using the condition of $ \frac{1}{4}$-symmetry $(p=2)$ one gets respectively:
$$Q_{\la[\a}Q_{\b]\a_1}=0,$$
$$Q_{\a_3[\a_2}Q_{\a_1][\a}Q_{\b]\la}+\omega_\a Q^{\rho}\,_{\a_3}Q_{\la [\rho}Q_{\b][\a_2}\omega_{\a_1]}-\omega_\b Q^{\rho}\,_{\a_3}Q_{\la [\rho}Q_{\a][\a_2}\omega_{\a_1]}=0.$$
Observe that since $Q_{\a\b}$ is totally orthogonal to $\text{span}\{\omega\}$, and $\omega$ is not null, the completely orthogonal parts with respect to $\omega$ of the two equations are
\bea
Q_{\la[\a}Q_{\b]\a_1}=0,\label{QQ}\\
Q_{\a_3[\a_2}Q_{\a_1][\a}Q_{\b]\la}=0.\nonumber
\eea
In fact, reasoning in a similar way, that is to say, contracting the $ \frac{1}{2p}$-symmetry condition 
$$  \n_{[\a_{2p}}\n_{\a_{2p-1}]}\ldots\n_{[\a_2}\n_{\a_1]}\left(\omega_{[\a}Q_{\b][\la}\omega_{\mu]}\right)=0$$
(after using the Ricci identity) with $\omega^{\mu}$ and $\omega^{\a_{2p}}$  and taking the completely orthogonal part with respect to $\omega$ one gets 
$$Q_{\a_{(2p-1)}[\a_{{(2p-2)}}}\ldots Q_{\a_1][\a}Q_{\b]\la}=0.$$
Any of these conditions leads to a $Q$ with matrix-rank 1. To see this, for {\em any} vector field $Y$ define the one-form $Q^{{Y}}$ by means of $Q^{(Y)}_\a:=Q_{\a\b}Y^\b$. Then, equation \eqref{QQ} leads to 
$$
Q^{(Y)} \wedge Q^{(\overline Y)} =0, \hspace{1cm} \forall Y,\overline Y \in TM.
$$
Hence, all the $Q^{(Y)}$ are proportional to each other for arbitrary vector field $Y$. This readily implies that  $Q_{\a\b}=\sigma_\a \rho_\b$ for some one forms $\sigma,\rho\in \Lambda(M)$. As $Q_{\a\b}$ is symmetric, necessarily $\sigma = A \rho$ for some scalar $A$, and $Q_{ab}=A\, \rho_\a \rho_\b$. Thus $R_{\a\b\la\mu}=4\omega_{[\a}\rho_{\b]}\rho_{[\la}\omega_{\mu]}$ and the manifold would be of type ${\rm II}$, contrary to the assumption. 
Similarly for larger values of $p$.
\begin{flushright}
  $\blacksquare$
\end{flushright}

Consider now the case of $\KS$ of type ${\rm I}_N$, so that $\omega$ is a lightlike one-form. 
%
Now the equations \eqref{rtb}-\eqref{ectb} reduce to
 \begin{eqnarray}
  R_{\a\b}&=&\omega_{\a} Q_{\b\rho}\omega^\rho+\omega_{\b} Q_{\a\rho}\omega^\rho-Q_\rho\,^\rho\omega_\a \omega_ \b;\label{rtb2}\\
S&=&2Q_{\rho\sigma}\omega^\rho\omega^\sigma.\label{ectb2}
 \end{eqnarray}

In the Lorentzian case, the existence of the null $\omega$ satisfying \eqref{pb} implies that $\omega$ is in fact a multiple {\em aligned null direction} (AND) of the Riemann tensor (a `RAND') so that the curvature tensor is algebraically special. For the definition and properties of AND, see \cite{OPP} and references therein. Furthermore, it is easy to prove from \eqref{pb} that the Weyl tensor satisfies
$$
\omega_\rho C^\rho\,_{\a[\b\la}\omega_{\mu]} =0
$$
so that $\omega$ is also a multiple AND for the Weyl tensor (a `WAND').

Instead of using the Weyl tensor, Kaigorodov \cite{KaiB} chose to provide a classification of type {\rm I}$_N$ based on the Ricci tensor. 
Observe that from \eqref{rtb2} one can easily deduce that the rank of the Ricci tensor matrix is less than or equal to 2. The classification combined this rank with the vanishing or not of the scalar curvature as follows. 
 \begin{enumerate}[(i)]
   \item when the Ricci tensor vanishes (rank 0) or its rank is equal to 1,
   \item when the rank of the Ricci tensor is equal to 2 and the scalar curvature does not vanish,
   \item  and finally when the rank is equal to 2 and the scalar curvature does vanish. 
 \end{enumerate} 
 For details consult \cite{KaiB}.

\subsubsection{Type {\rm II}}
 In this class of manifolds $\text{dim }\L=2$ and there exist two linearly independet one-forms $\omega,\rho \in \Lambda(M)$ that generate $\L$. 
%
    By Lemma \ref{co2}, the curvature tensor, the Ricci tensor and the scalar curvature assume the form:
 \begin{eqnarray}
 R_{\a\b\la\mu}&=& 4A\, \omega_{[\a}\rho_{\b]}\rho_{[\la}\omega_{\mu]};\label{ct}\\
 R_{\a\b}&=&A\left(2 (\omega_\nu \rho^\nu) \omega_{(\a}\rho_{\b)}-(\rho_{\sigma}\rho^{\sigma})\omega_{\a}\omega_{\b}-(\omega_{\rho}\omega^{\rho})\rho_{\a}\rho_{\b}\right);\label{rt}\\
 S&=&  2A\left[(\omega_\nu \rho^\nu)^2-\omega_{\rho}\omega^{\rho}\rho_{\sigma}\rho^{\sigma}\right].\label{ect}
 \end{eqnarray}
As a general result for type II we have the following.
\begin{theorem}[\cite{Kai,KaiA}]\label{teo2}
All $\KS$ spaces of type {\rm II} are semi-symmetric.
\end{theorem}
 
 {\bf Proof.}
Introducing the explicit expression for the curvature given in \eqref{ct} into the lefthand side of \eqref{RicR} the result is identically zero.
\begin{flushright}
{$\blacksquare$}\end{flushright}
Now take the expression \eqref{LDCai} 
and antisymmetrize in the indices $\la_{r-1}$ and $\la_r$. Since a $\KS$ space of type {\rm II} is semi-symmetric by Theorem \ref{teo2}, one gets that for type {\rm II}
\begin{equation}\label{ss}
  t_{\la_1\ldots[\la_{r-1}} \n_{\la_r]} R^\a\,_{\b\mu\nu}+t_{\la_1\ldots[\la_{r-1}\la_r]} R^\a\,_{\b\mu\nu}=0.
\end{equation}
\begin{coro}
    Let $(M,g)$ be a non-flat $\KS$ of type {\rm II} with $r\geq 2$. Then, if the $r$-covariant tensor  field $t^{(r)}$ in \eqref{LDC} is symmetric in its last two indices, it follows that
    $$t_{\la_1\ldots[\la_{r-1}} \n_{\la_r]} R^\a\,_{\b\mu\nu}=0$$
     and $t^{(r)}$ can be decomposed as:
    \be\label{decom}
    t^{(r)}=A_1\otimes \omega\otimes \omega  +   A_2\otimes (\omega \otimes \rho+\rho\otimes\omega) +  A_3 \otimes \rho\otimes\rho
    \ee
  for some $(r-2)$ covariant tensor fields $A_1, A_2$ and $A_3$.
    In particular, this is always true for type {\rm II} $KS_n$ spaces with $t^{(r-1)}=0$.
    \end{coro}

 {\bf Proof.}
If the r-covariant tensor $t^{(r)}$ in \eqref{LDC} is symmetric in its last two indices the first relation follows directly from \eqref{ss}. Also, contracting \eqref{gt} with $Z^{\la_r}\in\L^\perp$ and using \eqref{ct} one gets 
$$Z^{\la_r}t_{\la_1\ldots\la_{r-1}\la_r} =0, \, \hspace{1cm} \forall Z\in {\cal L}^\perp$$ 
since $(M,g)$ is non-flat. Given that $t^{(r)}$ is symmetric in its last two indices, the decomposittion \eqref{decom} follows at once. The last sentence is a consequence of Lemma \ref{p1}. The case where $t^{(r)} =0$ is included here with $A_1=A_2=A_3=0$.
\begin{flushright}
{$\blacksquare$}\end{flushright}


Type II $\KS$ semi-Riemannian manifolds can be naturally split into three different classes according to the level of degeneracy of ${\cal L}$. In particular we have --see also  \cite{KaiA}--
 \begin{itemize}
   \item Type {\rm II}$_\ep$ if $\L$ is non degenerate, so that it admits an orthonormal basis with $g(\omega,\omega):=\ep_\omega=\pm 1$,  $g(\rho,\rho):=\ep_\rho=\pm 1$ and $g(\omega,\rho)=0$. Here $\epsilon \in \{1,0,-1\}$ can be defined by $\ep:=(\ep_\omega +\ep_\rho)/2$.
   \item Type {\rm II}$_{N_\ep}$  if $\L$ is degenerate of rank 1,  i.e.,  if $\L$ admits a basis $\{\omega,\rho\}$ so that $g(\omega,\omega)= 0$,  $g(\rho,\rho):=\ep=\pm 1$ and $g(\omega,\rho)=0$.
   \item Type {\rm II}$_0$  if $\L$ is degenerate of rank 0 (the metric on $\L$ vanishes) so that $\L$ admits a basis of two null and mutually orthogonal null one-forms: $g(\omega,\omega)= 0$,  $g(\rho,\rho)=0$ and $g(\omega,\rho)=0$.
 \end{itemize} 
 Observe that in the Riemannian case only type {\rm II}$_1$ spaces are possible. In the Lorentzian case, types {\rm II}$_{-1}$, {\rm II}$_0$ and {\rm II}$_{N_{-1}}$ are not possible either. 
 

For type {\rm II$_\ep$ relations \eqref{rt} and \eqref{ect} simplify to 
 \begin{eqnarray}
 R_{\a\b}&=&\frac{S}{2} \left(\ep_\omega \omega_{\a}\omega_{\b}+\ep_\rho\rho_{\a}\rho_{\b}\right);\label{rt1}\\
 S&=&  -2\ep_\rho \ep_\omega A.\label{ect1}
 \end{eqnarray}
 The first of these states that the Ricci tensor is actually proportional to the metric on the 2-dimensional space $\L$. Furthermore, 
 \begin{propo} \cite{KaiA} \label{A1}
   Let $(M,g)$ be a non-flat $\KS$ of type {\rm II}$_\ep$. Then
   \begin{enumerate}[(a)]
     \item the scalar curvature is non-vanishing.
     \item The curvature tensor, the Ricci tensor and the scalar curvature satisfy the following equations:
     \begin{eqnarray}
        \frac{1}{2} S R_{\a\b\la\mu}=R_{\a\la}R_{\b\mu}-R_{\a\mu}R_{\b\la},\label{eqA11}\\
        R^\a\,_{\rho}R^{\rho}\,_{\b}=\frac{1}{2}S R^\a\,_\b,\label{eqA12}\\
         R_{\a\b}R^{\a\b}=\frac{1}{2}S^2,  \label{eqA13}\\
        R_{\a\b\rho\sigma}R^{\rho\sigma\la\mu}=SR_{\a\b\la\mu}, \label{eqA14}\\
        R_{\a\b\la\mu}R^{\a\b\la\mu}=S^2 \label{eqA15}.
     \end{eqnarray}
          \end{enumerate}
   \end{propo}
   {\bf Proof.}
   From \eqref{ct} and \eqref{ect1} one has
   $$
   R_{\a\b\la\mu}= 2\, \ep_\rho \ep_\omega S \, \omega_{[\a}\rho_{\b]}\rho_{[\la}\omega_{\mu]}
   $$
   so that $S=0$ would imply vanishing curvature. From this and \eqref{rt1} eqs. (\ref{eqA11}-\ref{eqA15}) follow.
\begin{flushright}
{$\blacksquare$}\end{flushright}

Consider now type {\rm II$_{N_\ep}$. The relations \eqref{rt} and \eqref{ect} simplify now to 
\begin{eqnarray}
 R_{\a\b}&=&-\ep A \omega_{\a}\omega_{\b};\label{rt2}\\
 S&=& 0. \label{ect2}
 \end{eqnarray}

\begin{propo} \cite{KaiA} \label{A2}
    Let $(M,g)$ be a $\KS$ of type {\rm II}$_{N_\ep}$. Then,
   \begin{enumerate}[(a)]
     \item The scalar curvature vanishes.
          \item The curvature tensor and the Ricci tensor satisfy
    $$ R_{\a\rho}R^\rho\,_\b=0 ;\,\,\,\, \,\,\,\, R_{\a\b}R^{\a\b}=0 ;\,\,\,\, \,\,\,\, R_{\a\b\la\mu}R^{\a\b\la\mu}=0.$$  
  \end{enumerate}
\end{propo}

{\bf Proof.}
  (a) is just formula \eqref{ect2}. (b) follows from \eqref{ct} and \eqref{rt2}.
     \begin{flushright}
{$\blacksquare$}\end{flushright}

For Lorentzian type {\rm II}$_{N_1}$, $\omega$ is a multiple AND \cite{OPP} of maximum multiplicity and the Riemann, Ricci and Weyl tensors are of the so-called type {\em N}. 
 
 Finally, consider type {\rm II}$_0$. Relations \eqref{rt} and \eqref{ect} become
 \be
 R_{\a\b}=0, \hspace{1cm} S=0.
 \ee
It is straightforward then to prove the following:
   \begin{propo} \cite{KaiA}\label{A3}
    Let $(M,g)$ be a $\KS$ of type {\rm II}$_0$. Then
   \begin{enumerate}[(a)]
     \item $(M,g)$ is Ricci flat.
       \item the curvature tensor coincides with the Weyl tensor and satisfies $R_{\a\b\la\rho}R^{\rho}\,_{\mu\nu\sigma}=0.$
     \end{enumerate}$\blacksquare$
 \end{propo}
%


\section{A model family: the class $\KS$ is not empty}\label{exp}
In this section, in order to prove that the homogeneous equation \eqref{LDC} has a non-empty set of solutions, an explicit family of semi-Riemannian metrics satisfying the general equations \eqref{LDC} and \eqref{nHLDC} for each $r\geq 1$ is presented in the following theorem:   
\begin{theorem}\label{exa}\cite{Jaen,SN2}
  Let $(\R^n,g)$ be an irreducible semi-Riemannian manifold of signature $(p+1,q+1)$, $p+q=n-2$ with Cartesian coordinates $\{u,v,x^2,\ldots, x^{n-1}\}$ endowed with the metric 
$$g=2du\left(dv+A_{ij}(u)x^ix^jdu\right)+\eta_{ij}dx^idx^j$$
where $\eta$ is the flat semi-Riemannian metric of signature $(p,q)$ on $\R^{n-2}$ and $A_{ij}(u)$ are functions of $u$ with $A_{ij}=A_{ji}$ ($i,j,\dots =2,\dots n-1$).
 
Choose fixed functions $B_{ij}(u)$ on $M$. Then, $(M,g)$ satisfies the linear differential condition on the curvature \eqref{nHLDC} with 
 $$
 t^{(m)}=a_{m}(u)\underbrace{du\otimes \ldots\otimes du}_{m}, \hspace{5mm} B=B_{ij}(u)\underbrace{du \otimes \dots \otimes du}_r \otimes \frac{\partial}{\partial v} \otimes dx^i\otimes dx^j  \otimes du$$
if and only if the functions $A_{ij}(u)$ satisfy 
\begin{equation}\label{eqej0}
    \frac{d\,^{r}A_{ij}}{du\,^r}+\sum_{m=1}^{r}a_{m}(u)\frac{d\,^{r-m}A_{ij}}{du\,^{r-m}}=B_{ij}(u)
\end{equation} 
for all $i,j$.
  
In particular, $(M,g)$ satisfies the homogeneous linear differential condition on the curvature \eqref{LDC}, and therefore is a $\KS$ space, if and only if all the functions $A_{ij}(u)$ satisfy the same and \underline{unique} homogenous ordinary differential equation 
  \begin{equation}\label{eqej}
    \frac{d^{r}y}{du^r}+\sum_{m=1}^{r}a_{m}(u)\frac{d^{r-m} y}{du^{r-m}}=0.
  \end{equation}
  \end{theorem}
{\bf Proof.}
If one takes the general co-basis for the cotagent bundle $T^*M$ given by:
$$\t^0=du;\,\,\,\,\,\,\,\, \t^1=dv+A_{ij}(u)x^ix^j du;\,\,\,\,\,\,\,\, \t^i=dx^i,$$
the metric can be written as
$$g=2\t^0\t^1+\sum_{i=2}^{n-1}\epsilon_i (\t^i)^2, \text{ with }(\epsilon_i)^2=1.$$
The connection 1-forms $\o^\a_\b$ can be computed from the formulae $\n g=0$, $d\t^\a=-\o^\a_\b\wedge \t^\b$ (see \cite{Car} for more information) and the fact that  $\t^0=du$ is a null parallel one-form. Then, the only non-vanishing connection 1-forms are:
\begin{eqnarray*}
 \o^1_i=A_{ik} x^k\t^o=-\o^i_0.
 \end{eqnarray*}
 Calculating the curvature 2-forms $\Omega^\a_\b$ by the formulae $\Omega^\a_\b=d\o^\a_\b+\o^\a_\rho\wedge \o^\rho_\b$, one deduces that the only non-vanishing  components of the curvature $R$ are:
 $$R^1\,_{ij0}=A_{ij}(u).$$
 It is easy to show that, therefore, the non-vanishing components of $\n R$, $\n^2 R$, \ldots, $\n^r R$ in this co-basis are, respectively:
 \begin{eqnarray*}
   \n_0 R^1\,_{ij0}&=&\frac{d}{du}{A}_{ij}(u);\\
    \n_0 \n_0 R^1\,_{ij0}&=&\frac{d^2}{du^2}{A}_{ij}(u);\\
   & \vdots&\\
    \underbrace{ \n_0 \ldots \n_0}_r R^1\,_{ij0}&=&\frac{d\,^r}{du\,^r}{A}_{ij}(u).
 \end{eqnarray*}
 In conclusion, the only non-vanishing components that survive for the terms $\nabla^{m} R$ on the left-hand side of the homogeneous lineal differential conditions \eqref{nHLDC} in this co-basis are those associated to 
 $$\underbrace{\t^0\otimes \ldots\otimes \t^0}_m \otimes \partial_v \otimes \t^i \otimes \t^j \otimes \t^0$$
 and therefore the only way that \eqref{nHLDC} can be satisfied is that the tensor fields $t^{(m)}$ have the form 
 $$
t^{(m)}=a_{m}(u)\underbrace{du\otimes \ldots\otimes du}_m .$$
Then, \eqref{nHLDC} reduces in this co-basis to \eqref{eqej0}.
%
If the fixed tensor field $B=0$, and as all the $t^{(m)}$ do not depend on $i,j$, then all the $A_{ij}$ must satisfy the unique ODE \eqref{eqej}. 

 \definecolor{shadecolor}{rgb}{1,0.8,0.3} 
\begin{flushright}
{$\blacksquare$}\end{flushright}
 This family contains an example of a solution for a specific homogeneous linear differential equation of order $r\geq 1$, for each $r$. Observe that each $t^{(m)}=a_m(u)\t^0\otimes\ldots\otimes\t^0$ is a fully symmetric recurrent tensor field $\forall m=1,\ldots,r$, since $\n t^{(m)}=\displaystyle\frac{\dot{a}_m}{a_m}\t^0\otimes t^{(m)}$, where $\dot{\,\,\, }$ means derivation with respect to the variable $u$. This raises the questions of whether or not, in general, the tensor fields $t^{(m)}$ will have any type of symmetry or recurrent properties.
 
 All the manifolds in the new family given in Theorem \ref{exa} possess a parallel lightlike one-form $du$, whose contravariant version is $\partial_v$, so that they realize all signatures {\em except the Riemannian one}. This model family in Theorem \ref{exa} has vanishing scalar curvature $S=0$, and the Ricci tensor field is given by 
 $$Ric=\eta^{ij} A_{ij}(u)du\otimes du$$
 where $\eta^{ij}$ is the contravariant, inverse, metric of $\eta_{ij}$. In the classification of the previous section, they belong to one of the following cases, always with $\omega =\t^0 =du$:
 \begin{itemize}
 \item type {\rm I}$_N$ with the additional property that $Q_{\a\b}\omega^\b =0$, where $Q_{\a\b}$ corresponds essentially to $A_{ij}$.
 \item type {\rm II}$_{N_1}$ if rank$(A_{ij})=1$.
 \end{itemize}

As a final comment, it is remarkable that there exists a one-to-one correspondence between the $\KS$ model family described in Theorem \ref{exa} and the family of order-$r$ linear ordinary differential equations.
 
 \section{Particular cases of special relevance}\label{sec:cases}
 In this section a compendium of the most relevant particular families of semi-Riemannian manifolds contained in the $\KS$ manifolds that have been already studied in the literature is given, with particular emphasis in their (non)-existence and their known properties.
For more than 100 years a vast work has been developed characterizing and classifying semi-Riemannian manifolds with conditions on the curvature. The outstanding example is the class of {\em locally symmetric spaces}, defined by the condition 
$$\n R=0$$
so that $R$ is parallel propagated along any direction.
Introduced by Shirokov in \cite{Sh} these were thoroughly studied and classified by Cartan in 1926-27 \cite{CA,CA2} in the Riemannian case, and later in 1970 by Cahen and Wallach in the Lorentzian case \cite{CW}. 

The next important family is the set of {\em recurrent spaces}, defined by the condition 
$$\n R=\sigma\otimes R$$ 
for some 1-form $\sigma\in \Lambda(M)$, so that $\n R$ is being propagated parallel to itself along any direction.
First introduced by Ruse in 1946 \cite{Ru}, they contain the locally symmetric ones and thus they are referred to as {\em proper recurrent space} if $\o\neq 0$. Locally symmetric and proper recurrent $(M,g)$ belong to the family of $\KS$ spaces for $r=1$, the former with $t^{(1)}=t^{(r)}= 0$.

More generally, {\em $r^{th}$-order recurrent spaces}, or shortly {\em $r$-recurrent spaces},  were introduced in \cite{CS}, and are defined by the condition 
$$\n^{r}R=-t^{(r)} \otimes R$$ 
where $t^{(r)}$ is a $r$-covariant tensor field. If  $t^{(r)}\neq 0$, they are called {\em proper $r^{th}$-recurrent spaces} while if $t^{(r)}=0$, so that the curvature condition reads simply
$$\n^{r}R=0,$$
they are called {\em $r^{th}$-order symmetric spaces} --or shortly {\em $r$-symmetric spaces}-- and were introduced\footnote{Kaigorodov already defined these spaces in \cite{Kai} but with the added condition that the subspace $\mathcal{L}$ of Definition \ref{Lspace} is non-empty. Examples of $r$-symmetric spaces were given in \cite{Kai} and in \cite{E-K}.} in \cite{SN1}.  A $r^{th}$-order symmetric space is {\em proper }if it satisfies $\n^{r}R=0$ but $\n^{r-1}R\neq 0$. Proper $r$-recurrent and $r$-symmetric spaces belong to the class $\KS$ with $t^{(m)}=0$ for all $m=1,\dots, r-1$ in the former case, and also with $t^{(r)}=0$ in the latter.

As a corollary of Lemma \ref{tr=0}, one can state the following: 
\begin{coro}[\cite{SN1,TA}]\label{rrrs}
At generic points, $r^{th}$-order recurrent semi-Riemannian manifolds are $r^{th}$-order symmetric.
\end{coro}


Important known results about locally symmetric, recurrent, proper $r$-symmetric and proper $r$-recurrent spaces are summarized in what follows. Starting with Riemannian Geometry, 
traditional results show that there are no proper $r$-symmetric Riemannian spaces with $r>1$ \cite{LI2,NO,TA}; and letting aside obvious flat extensions of 2-dimensional recurrent spaces, for $n>2$ there are no proper recurrent Riemannian manifolds \cite{moG,WAL}. Furthermore, there are no proper $2^{nd}$-order recurrent Riemannian spaces, all of them being locally symmetric. Summarizing:
\begin{theorem}\label{RieCase}
  Let $(M,g)$ be a Riemannian manifold of dimension $n$. Then:
  \begin{enumerate}[(a)]
    \item \cite{LI2,TA} if it is $r^{th}$-symmetric, or
    \item \cite{Ro2,Ro3}  if it is recurrent, or proper $2^{nd}$-order recurrent and $n>2$,
  \end{enumerate}
  then it is actually locally symmetric.
\end{theorem}
Concerning recurrent spaces for signatures other than Riemannian one has the following theorem due to Walker \cite{WAL} that provides a classification of recurrent semi-Riemannian manifolds: 
\begin{theorem}\cite{WAL}\label{re}
  Letting aside obvious flat extensions, for $n>2$ any proper recurrent semi-Riemannian manifold of signature (p+1,q+1) satisfying $\n R=\sigma \otimes R$ with $\sigma\neq 0$:
  \begin{enumerate}[(a)]
\item is not reducible;
\item has a parallel lightlike one-form $\o$ ($\n \o=0$ and $g(\o,\o)=0$) such that $\sigma\wedge \omega=0$;
\item belongs to the model family of $\KS$ given in Theorem \ref{exa} with $\o =du$, $A_{ij}(u) = F(u) M_{ij}$, $M_{ij}\in \mathbb{R}$ with rank$(M_{ij})\geq 1$ and $F(u)$ is a non-constant, otherwise arbitrary, function of $u$. The recurrence one-form is $\sigma=-\ds\frac{\dot{F}(u)}{F(u)}du$ with $\dot{F}(u):=\ds\frac{dF}{du}$. 
\end{enumerate}  
\end{theorem}
 Observe that this theorem shows that there are no proper recurrent Riemannian manifolds, in accordance with Theorem \ref{RieCase}. 
 
 Notice that from Lemma \ref{p1}, for a proper $r^{th}$-recurrent space ($r\geq 2$) $t^{(r)}$ is always symmetric in its last two indices.
 \begin{propo} \cite{Kai}\label{p2}
   If  $(M,g)$ is proper $r^{th}$-order recurrent, including $r^{th}$-order symmetric, with $r\geq 2$ then it is $\frac{1}{2p}$-symmetric for $p=r/2$ if $r$ is even, and $p=(r+1)/2$ if it is odd.
 \end{propo} 
 {\bf Proof.}
 If the manifold is $r^{th}$-symmetric the proof follows immediately for $p=r/2$ if $r$ is even, and $p=(r+1)/2$ if it is odd.  Suppose then that the manifold is $r^{th}$-recurrent. If $r$ is an even number, take $p=r/2$. The condition for $\frac{1}{2p}$-symmetry given is then
   $$ \n_{[\la_1}\n_{\la_2]}\ldots \n_{[\la_{r-1}}\n_{\la_{r}]} R^\a\, _{\b\la\mu}=0 $$
and as the manifold is $r$-recurrent this equation is equivalent to
$$ t_{{[\la_1}{\la_2]\ldots[\la_{r-1}}{\la_{r}]}} R_{\a\b\la\mu}=0, $$
which is trivially satisfied for $r^{th}$-recurrent manifolds by Lemma \ref{p1}. If $r$ is odd, then take $p=\frac{r+1}{2}$ to obtain that the condition for $\frac{1}{2p}$-symmetry reads
   $$ \n_{[\la_1}\n_{\la_2]}\ldots \n_{[\la_{r}}\n_{\la_{r+1}]} R_{\a\b\la\mu}=0$$
which is rewritten by $r^{th}$-recurrence as
$$ \n_{[\la_1}t_{\la_2]\ldots[\la_{r}\la_{r+1}]} R_{\a\b\la\mu}=0, $$
which is again trivially satisfied due to Lemma \ref{p1}.

\begin{flushright}
  $\blacksquare$
\end{flushright}

\begin{coro} \label{c3}
  There are neither $r^{th}$-order recurrent nor r-symmetric semi-Riemannian manifolds of type ${\rm I}_\ep$. \end{coro}
  {\bf Proof.}
This follows from Proposition \ref{p2} and Theorem \ref{p5}.
\begin{flushright}
  $\blacksquare$
\end{flushright}

 There is no formal classification of $r^{th}$-order recurrent spaces with $r\geq 2$.  Thompson \cite{To4} proved that if a manifold is proper $2^{nd}$-order recurrent and conformally flat (which implies that the scalar curvature is zero due to Chowdhury \cite{Ch}), then it is recurrent. In \cite{To3}, he proved that if the scalar curvature is not zero, and a proper $2^{nd}$-recurrent manifold is either 3-dimensional, 4-dimensional and Lorentzian, or its Ricci tensor is definite, then it is recurrent. Also, he proved in \cite{To2} and \cite{To1} that a 4-dimensional proper $2^{nd}$-order recurrent Lorentzian manifold is complex-recurrent. With all this information at hand, in 1972 together with McLenaghan they proved that:
 
 \begin{propo} \cite{To5}\label{prop1}
   A proper $2^{nd}$-order recurrent Lorentzian manifold of dimension 4 is non-corformally flat  with vanishing scalar curvature, and there exist coordinates $\{u,v,x,y\}$ such that the metric takes the form:
 $$g=2du\left( dv+H(u,x,y)du\right)+dx^2+dy^2,$$
 where $\p_v$ is a parallel lightlike vector field, 
 $$H(u,x,y)=h(u)(x^2+y^2)+2f(u)\left[cos(\b(u))(x^2-y^2)-2\sin(\b(u))xy\right],$$
 with $f(u)$ an arbitrary function and $\b(u)$ and $h(u)$ satisfying the equations
 $$\begin{cases}
 \dot\b=\kappa/f^2, &\kappa\in\R .\\
 \ddot h -h/f(\ddot f-\kappa^2/f^3)=0.
 \end{cases}
 $$
 The manifold is Ricci-flat if and only if $h=0$, and it is recurrent if and only if $\dot \b(u)=\kappa =0$ and $h=\lambda f$ for some $\lambda\in\R$.
 \end{propo}  
Again, these spaces belong to the model family of Theorem \ref{exa}. In this case , the relation between the metric in Proposition \ref{prop1} and that of Theorem \ref{exa} is given by
\begin{align*}
 A_{xx}(u)&=h(u)+2f(u)\cos(\b(u)),\\
 A_{yy}(u)&=h(u)-2f(u)\cos(\b(u)),\\
 A_{xy}(u)&=A_{yx}(u)=-4f(u)\sin(\b(u)).
\end{align*}
Using the equations in Proposition \ref{prop1} that relate $f(u)$, $h(u)$ and $\b(u)$ it can be proven that the above functions $A_{ij}(u)$ do satisfy the conditions given in Theorem \ref{exa} for $n=4$ and $(p+1,q+1)=(1,3)$, with 
$$t^{(2)}=\ds\left(\frac{\ddot{f}}{f}-\frac{\kappa^2}{f^4}\right)du\otimes du.$$
In the solution given in Proposition \ref{prop1} there are an arbitrary funtion $f(u)$ and four arbitrary constants, $\kappa$ and the ones obtained after solving the equations for the functions $\beta$ and $h$. However, $f(u)$ and $\kappa$ define the tensor field $t^{(2)}$ in the equation of second recurrence, and then only the sign of $\kappa$ matters because, if $\kappa \neq 0$, by taking $\bar{f}(u)=\frac{f}{\sqrt{|\kappa|}}$, the constant $\kappa$ becomes equal to $\pm 1$. Thus, there are 3 essential parameters plus one sign in the solution. On the other hand, in the model family of solutions of Theorem \ref{exa}, since $n=4$ and $r=2$ there would be a priori 6 arbitrary parameters in the metric. This apparent conflict can be easily resolved by performing an appropriate change of coordinates that preserves the form of the metric in Theorem \ref{exa} ---see for instance Claim 5.20 in \cite{BSS3} or \cite{Blau}--- which allows to remove two of these constants.

Passing now to the cases of $r$-symmetry for $r\geq 1$, the full classification of the Lorentzian symmetric spaces is known.
\begin{theorem}\cite{CW} 
Any simply-connected Lorentzian symmetric space $(M,g)$ is isometric to the product of a simply-connected Riemannian symmetric space and one of the following Lorentzian manifolds:
\begin{enumerate}[(a)]
 \item $(\mathbb{R},-dt^2)$
\item the universal cover of d-dimensional de Sitter or anti-de Sitter spaces, $d\geq 2$, 
\item a metric of $\KS$ type in the model family of Theorem \ref{exa} with all $A_{ij}$ constant and $\eta_{ij}=\delta_{ij}$, the Kronecker delta.
\end{enumerate}
The $d$-dimensional Lorentizan manifolds with the metric of (c) are sometimes called Cahen-Wallach spaces and denoted by $CW^d$ (observe that $CW^2 =\mathbb{L}^2$ is just the two-dimensional Minkowski space, as $A_{ij}$ necessarily vanish). 
Therefore, if a Lorentzian symmetric space admits a parallel lightlike vector field, then it is locally isometric to the product of a $d$-dimensional Cahen-Wallach space and an $(n-d)$-dimensional Riemannian symmetric space with $d\geq 2$.
\end{theorem}
Locally symmetric spaces of signatures other than positive definite and Lorentzian are not so well understood. Important advances can be found in \cite{CP,CK,KO,KO1} and references therein.

In the case of $r^{th}$-order symmetric spaces some progress has been made in this century for $r=2,3$ in the Lorentzian case. For $r=2$, in \cite{SN1} it was proven that such spaces admit a lightlike parallel vector field, and latter in \cite{BSS2,BSS3} the global classification was given. The local classification was also given in \cite{AG} using Lorentzian holonomy techniques \cite{G2}. Some years later $3^{rd}$-order Lorentzian symmetric spaces were also locally classified in \cite{G3} using the same techniques as in \cite{AG}. The following theorem combines all these results
\begin{theorem} \cite{BSS3,AG,G3} \label{sesi}
  Letting aside obvious locally-symmetric extensions, for $n>2$ any proper $2^{nd}$-order or $3^{rd}$-order symmetric Lorentzian manifold:
  \begin{enumerate}[(a)]
    \item is not reducible;
    \item possesses a parallel lightlike one-form $\o$: $\n \o=0$;
    \item belongs to the $\KS$ model family of Theorem \ref{exa} with $\o =du$, $\eta_{ij}=\delta_{ij}$,
    $$
    A_{ij}(u) = D_{ij}u^2+B_{ij}u+C_{ij}
    $$
    with $D_{ij},B_{ij},C_{ij} \in \mathbb{R}$ and rank$(D_{ij})\geq 1$ for 3-symmetry, while $D_{ij}=0$ and rank$(B_{ij})\geq 1$ for 2-symmetry.
    \end{enumerate}
  \end{theorem}
  
The global result is
\begin{theorem} \cite{BSS3,G3}\label{global}
Any geodesically complete and simply connected proper $2^{nd}$- or $3^{rd}$-order symmetric Lorentzian manifold is globally isometric to a direct product of a non-flat Riemannian symmetric space and the 2- or 3-symmetric Lorentizan manifolds of Theorem \ref{sesi}.
\end{theorem}
 
It is clear that the model families in Theorem \ref{exa} with $A_{ij}(u)$ polynomials of order $r-1$ are $r$-symmetric for any $r\in \mathbb{N}$. However, there is no proof that these exhaust the entire family for $r\geq 4$. According to \cite{G3}, Lorentzian holonomy techniques \cite{G2} are no longer usable for this purpose if $r\geq 4$, hence, the main line to be pursued is to try and prove that,  for any $r$, there is a parallel lightlike one-form. 
%

A question that springs to mind is whether or not a $\KS$ space can satisfy more than one relation of type \eqref{LDC}. This can be answered in the affirmative by the following result.
\begin{coro}\label{cor2}
Letting aside obvious flat extensions, a Lorentzian manifold $(M,g)$ of dimension $n>2$ that satisfies both the  recurrent condition and the $k^{th}$-order symmetric condition on the curvature for $k=2,3$ is given by the solutions of Theorem \ref{re} with $F(u)$ a polynomial of degree $k-1$ with leading coefficient $a_{k-1}=1$.
\end{coro}
{\bf Proof.}
If $(M,g)$ is Lorentzian and recurrent, it has a parallel one-form $\o$ (Theorem \ref{re}). If it is also $k$-symmetric, $k\in \{2,3\}$, it also possesses a parallel one-form $\tilde\o$ (Theorem \ref{sesi}). Hence, the first thing one needs to prove is that these two parallel directions coincide. Suppose, on the contrary,  that they are linearly independent. Then, they would span a 2-dimensional subspace of parallel one-forms, with signature necessarily Lorentzian, implying the existence of a parallel timelike one-form $U$ (and also a parallel spacelike one-form). This would imply that $(M,g)$ is decomposable as $(M_1\times M_2, g_1\oplus g_2)$, with  $g_1=\ds\frac{1}{g(U,U)} U\otimes U $ and $g_2$ a Riemannian metric. Since $(M,g)$ is both recurrent and $k^{th}$-symmetric for $k=2,3$, so would be $(M_2,g_2)$, and since $(M_2, g_2)$ is Riemannian, $(M,g)$ would be locally symmetric by Theorem \ref{RieCase}. In conclusion, if a Lorentzian manifold is proper recurrent and proper $k^{th}$-symmetric for $k=2,3$, the parallel one-forms $\o$ and $\tilde\o$ must be proportional --and can be chosen to be the same. Once this is established, the result follows at once from Theorems \ref{re} and \ref{sesi}.
\begin{flushright}
  $\blacksquare$
\end{flushright}

For general $r^{th}$-order symmetric semi-Riemannian manifolds one can prove the following interesting intermediate result 
\begin{propo}\label{r-sym}
Let $(M,g)$ be an irreducible $r^{th}$-order symmetric non-flat semi-Riemannian manifold.  Then either
\begin{enumerate}[(i)]
\item there exist a parallel lightlike vector field and $g\left(\nabla^{r-1} R, \n^{r-1} R \right)=0$, or
\item $g\left(\nabla^{r-2} R, \n^{r-2} R \right)=${\rm constant}.
\end{enumerate}
\end{propo}
{\bf Proof.}
Applying Lemma \ref{nnT} in Appendix \ref{nnT=0} to $T=R$, one sees that case $(i)$ in that Lemma is not possible due to Proposition \ref{hv4}. The result then follows immediately from cases $(ii)$ and $(iii)$ in Lemma \ref{nnT}.
\begin{flushright}
  $\blacksquare$
\end{flushright}
The importance of this result lies in the fact that, if case (i) in Proposition \ref{r-sym} holds, then a Lorentzian $r$-symmetric $(M,g)$ will be a Brinkmann manifold \cite{BR,BSS3}, and then using calculations such as in \cite{BSS2,BSS3,AG,G3} one may try to prove that such manifolds belong to the model family of Theorem \ref{exa}.

Strong results can also be obtained if there are generic points (Appendix \ref{sec32}) in the manifold. For instance
\begin{theorem} \cite{SN1} \label{gp}
  At generic points, all semi-symmetric semi-Riemannian manifolds are of constant curvature.
\end{theorem}
And taking into account the results in Appendix \ref{nnT=0} one obtains
\begin{propo}\label{kls}
  If a $r^{th}$-order recurrent  or a $r^{th}$-order symmetric semi-Riemmanian manifold $(M,g)$ has a generic point, then it is locally symmetric. In other words, a (proper) $r^{th}$-symmetric space has no generic points.
\end{propo}
This proposition follows from Corollary \ref{rrrs} and Theorem \ref{rsls} by  taking $T=R$ in Theorem \ref{rsls}.
Therefore, since locally symmetric spaces are semi-symmetric, from Propositions \ref{kls} and Theorem \ref{gp} one concludes that:
\begin{theorem}\label{kls2}
  At generic points, all $r^{th}$-order recurrent and all $r^{th}$-order symmetric semi-Riemannian manifolds are actually of constant curvature.
\end{theorem}

The full classification of semi-symmetric 4-dimensional Lorentzian manifolds can be found in \cite{ES,Aman}.

\section{Relevance in Physics}\label{sec5}
The physics of the gravitational field is described on Lorentzian manifolds, also called {\em spacetimes}. Actually, the gravitational field itself corresponds to the curvature of the spacetime. Therefore, results concerning $\KS$ Lorentzian manifolds have undoubtedly relevance for many branches of gravitational physics. A long list of applications was detailed in the Introduction of \cite{SN1}. They were listed for $2^{nd}$-order symmetric spacetimes, but they remain valid, {\em mutatis mutandis}, for the general class of $\KS$ Lorentzian manifolds. Thus, $\KS$ spacetimes may have some applications in (see \cite{SN1} and references therein)
\begin{itemize}
\item building conserved quantities of the gravitational field depending on higher-order derivatives of the  curvature, via for instance the so-called `super-energy' tensors \cite{S};
\item simplifying the expansions in Riemann normal coordinates, if these are useful;
\item the regularization of quantum fluctuations via curvature counter-terms;
\item the study of `Penrose limits';
\item as solutions of higher-order Lagrangian theories, in particular in some supergravity theories, string theory, and their relatives. Gauss-Bonnet gravity theories are of particular relevance in this case.
\end{itemize}
I would like to comment further in the last two points. Let me start with the case of Gauss-Bonnet gravity, see the review \cite{FCCM} and references therein. These are theories where the Lagrangian density contains not only the Einstein-Hilbert term $S$, but also an additional term proportional to the Gauss-Bonnet invariant $G$ defined in \eqref{e3}. This additional term is relevant for dimension $n>4$, as it provides field equations different from the Einstein equations but still of second-order and with two degrees of freedom. If the dimension is $n=4$, the integral of $G$ is a topological invariant —it equals the Euler characteristic of the manifold---and thus its contribution to the field equations vanishes upon extremization. As shown in Corollary \ref{coro2} $(a)$, the Gauss-Bonnet scalar vanishes for all $\KS$ spaces in any dimension, except possibly for type 0, if this exception happens to exist at all. Therefore, the $\KS$ spaces seem to be solutions of any Gauss-Bonnet extension of General Relativity in arbitrary dimension.

With regard the other point, let me first of all recall that the spacetimes in the model family of Theorem \ref{exa}, in the Lorentizan case (so that $\eta_{ij} =\delta_{ij}$) are called {\em plane-wave spacetimes}, or simply plane waves, in the physics literature \cite{EK,Blau}. They are geodesically complete \cite{CFS} for good-behaved $A_{ij}(u)$.
\begin{theorem}[Every Lorentzian $(M,g)$ has a ``Penrose limit'']\label{Plimit}\cite{Pen}
To every Lorentzian $(M,g)$ and choice of a null geodesic $\gamma$ one can associate a unique plane-wave space-time (i.e., one of the models in Theorem \ref{exa} ) called its ``Penrose limit'' on $\gamma$. 
\end{theorem}
The original construction \cite{Pen} of the Penrose limit amounts to selecting appropriate local coordinates adapted to the chosen null geodesic, performing a re-scaling with a parameter $\lambda$ and then taking the limit $\lambda \rightarrow 0$. This is why they are called `limits' in the first place. Penrose limits retain the information about geodesic deviation of the original spacetime along the chosen null geodesic $\gamma$, and nothing more. However, they can be determined in a covariant way without taking any limit, see \cite{Blau} and references therein: just compute the geodesic deviation (in orthogonal directions) along the selected null geodesic in the original space-time and this is identical to the null geodesic deviation in the corresponding Penrose limit plane wave metric. Such deviation is characterized by the symmetric $A_{ij}(u)$ in Theorem \ref{exa}, where $u$ plays the role of parameter along the null geodesic on the original spacetime. \cite{Blau}

This shows the information about the original metric that the Penrose limit encodes is precisely that of geodesic deviation along the selected null geodesic: in other words, it gives the curvature properties around the null geodesic to first non-trivial order ---in analogy with what Riemannian normal coordinates do at any given point \cite{BFW}. Furthermore, 
when approaching any spacetime singularity --the edge of incomplete geodesics \cite{SG,BEE}-- the Penrose limits acquire some scale-invariant properties and have a universal power-law behaviour related to the growth of the curvature in the original $(M,g)$ \cite{Blau}.

An obvious question is what properties of the original spacetime, if any, are preserved in its plane-wave Penrose limit. The study of the properties preserved if one considers a one-parameter family of space-times and takes some limit on the defining parameter was curried out long ago by Geroch \cite{Ge}. He defined a property shared by the entire one-parameter family of spacetimes {\em hereditary} if all the limits of this family also have this property. Following this idea, Blau \cite{Blau} states that a property of a Lorentzian manifold $(M,g)$ is hereditary if {\em all} its Penrose limits also have this property.

Any tensor field concomitant of the curvature --that is, any tensor field constructed from $R$, its derivatives, the metric and its inverse-- that vanishes is hereditary \cite{Ge}. For instance, every Penrose limit of de Sitter and anti-de Sitter spacetimes, for arbitrary null geodesics, is a flat space-time. One can conjecture that, conversely, if a space-time has all of its Penrose limits flat, then it must be of constant curvature. In general \underline{trace} properties of the Riemann tensor are lost in the limits. For instance, Einstein spaces have Ricci-flat Penrose limits --as all plane waves have vanishing scalar curvature. Using Geroch's ideas it seems somehow obvious to claim
\begin{quote}
{\it Every $\KS$ proper of order $r$ Lorentzian $(M,g)$ has as its {\em universal} Penrose limits the plane-wave models of Theorem \ref{exa} satisfying an appropriate homogeneous ODE \eqref{eqej} of order $r$---including a flat space-time as a particular sub-possibility }.
\end{quote}
The question is now to ascertain how much ``room'' is left between the plane-wave models of Theorem \ref{exa}, which are the Penrose limits, and the general Lorentzian manifolds of type $\KS$.

\section{Discussion: Open questions and conjectures}\label{sec6}
From the long list of results herein presented there are some facts that catch one's eye. The two most obvious, and probably more relevant, ones are 
\begin{itemize}
\item the absence of non-locally-symmetric Riemannian $\KS$ manifolds 
\item the existence of a parallel lightlike vector field in all known examples, in particular in the model family of Theorem \ref{exa}
\end{itemize}

This leads to the first open question
\begin{oq}[Riemannian manifolds in the new family?]
Are there any $\KS$ Riemannian manifolds other than the locally symmetric ones?
\end{oq}
All known partial results answer this question in the negative. Thus, Theorem \ref{RieCase} forbids the existence of proper recurrent, proper $r^{th}$-symmetric, and proper 2nd-order recurrent Riemannian spaces.  Furthermore, Proposition \ref{NoEinstein} ensures that no Riemannian $\KS$ Einstein space with non-empty $\L$ exists. The existence of generic points also leads to the absence of proper $r$-recurrent $\KS$ due to Theorem \ref{kls2}. And, of course, the model family of Theorem \ref{exa} does include all signatures except the Riemannian one. This leads to my first conjecture.
\begin{conj}\label{nelv}
There are no Riemannian $\KS$ manifolds of any order $r$, apart from the locally symmetric spaces.
\end{conj}
Notice that the combination of Theorem \ref{red} and Theorem \ref{RieCase} implies that, in dealing with this conjecture, only irreducible manifolds must be studied due to the de Rham decomposition \cite{RH}. 

On the other hand, all known (non-locally symmetric) examples of proper $\KS$ manifolds belong to the model family of Theorem \ref{exa}, which raises the second open question
\begin{oq}[Do the models in Theorem \ref{exa} exhaust the $\KS$ spaces?]\label{oq2}
Are there any $\KS$ semi-Riemannian manifolds other than the explicit models of Theorem \ref{exa} and the locally symmetric ones?
\end{oq}
Note that the explicit model family in Theorem \ref{exa} has a parallel null vector field. Thus, a second conjecture emerges.
\begin{conj}\label{constantX} 
  If $(M,g)$ is a non-locally symmetric $\KS$ semi-Riemannian manifold then it admits a parallel lightlike vector field.
\end{conj}
Again, all known partial results support this conjecture. Generally the model family in Theorem \ref{exa}, but also Theorem \ref{re} for the recurrent case, Theorem \ref{sesi} for $2^{nd}$- and $3^{rd}$-order symmetry, as well as Proposition \ref{prop1} for $2^{nd}$-order recurrent Lorentzian 4-dimensional manifolds and the partial result in Proposition \ref{r-sym} 
for general $r^{th}$-order symmetry.

It would be very important to settle Conjecture \ref{constantX}. In the Lorentzian case, and if the conjecture holds, then one can restrict the study to the Brinkmann metrics \cite{BR,BSS3} and the analysis would be much simplified. Furthermore, if Conjecture \ref{constantX} is resolved in the affirmative, then Conjecture \ref{nelv} will be immediately true. 

The validity of Conjecture \ref{constantX} would probably lead to a negative answer for the Open Question \ref{oq2}, and in that case one can formulate the more ambitious speculation.
\begin{conj}
All $\KS$ semi-Riemannian manifolds are contained in the family of Theorem \ref{exa}. In particular, for the Lorentzian case, they coincide with their own Penrose limits built with respect to the null geodesics tangent to the parallel vector field.
\end{conj} 

If anyone wishes to attack the previous open problems by steps, the simplest unresolved case is that of proper $2^{nd}$-order recurrent semi-Riemannian spaces. All partial known results are in agreement with the conjecture, as detailed next
\begin{propo}[2nd-order recurrent case: known results]\cite{To2,To1,To5,CS}
Let $(M,g)$ be a proper 2nd-order recurrent $n$-dimensional semi-Riemannian manifold, then
\begin{enumerate}
\item if it is also conformally flat, it is actually recurrent
\item if the scalar curvature does not vanish and 
\begin{itemize}
\item $R_{\mu\nu}$ is definite, then it is actually recurrent 
\item $g$ is Lorentzian, then  it is actually recurrent 
\end{itemize}
\item if $n=4$ and $g$ is Lorentzian, then the only solutions are the explicit plane wave models shown in Theorem \ref{exa}
\end{enumerate}
\end{propo}
Recall that the recurrent case is fully resolved and the solutions, given in Theorem \ref{re}, do belong to the model family of Theorem \ref{exa}. All in all, it seems that this case might not be too difficult to settle. Another important partial unresolved problem is that of Lorentizan $r$-symmetry for $r\geq 4$.

%
%
%
%
%

 

\section*{Acknowledgements}
I am grateful to S. Galaev for bringing ref.\cite{Kai} to my attention and for comments on the manuscript, to M. S\'anchez for many suggestions and helpful comments, and to M. Blau for some clarifications concerning Penrose limits. Work supported by Basque Government grant IT1628-22, and by Grant PID2021-123226NB-I00 funded by the Spanish MCIN/AEI/10.13039/501100011033 together with “ERDF A way of making Europe”.

\section*{Appendices}\label{classi}
In this section classical methods related to this type of problems are explained and known results presented. These methods could be helpful in the study of the properties of the new family of $\KS$ semi-Riemannian manifolds --those satisfying \eqref{LDC}.

\appendix

\section{Generic Points}\label{sec32}
Let $(M,g)$ be a semi-Riemannian manifold. Consider the curvature endomorphism on the space of two-forms (see, for example, \cite{ON}):
\bean {\cal R}: &\Lambda^2 M &\longrightarrow \Lambda^2 M\\
&\Omega &\mapsto R(\Omega)
\eean
defined by
$$
R(\Omega)_{\la\mu} := {R^{\a\b}}_{\la\mu}\Omega_{\a\b}.
$$
\begin{defn}
 A point $x\in M$ in a semi-Riemannian manifold $(M,g)$ is {\rm generic} if the curvature endomorphism at $x$ is  an isomorphism, that is to say, if ${\cal R}|_{x}$ is non-singular and has an inverse.
\end{defn}
Let $x\in M$ be a generic point. Then, there exists the inverse endomorphism ${\cal R}_{|x}^{-1}$ of ${\cal R}|_{x}$ --denoted in abstract index form by ${(R_{|x}^{-1})^{\a\b}}_{\la\mu}$-- at the point $x$. In fact, by the smoothness of $R$, the inverse  ${\cal R}^{-1}$ is defined in some neighbourhood $U_x$ of $x$. This endomorphism satisfies: 
\begin{equation}\label{RR-1}
  {R_{|y}~^{\a\b}}_{\rho\sigma}{(R_{|y}^{-1})^{\rho\sigma}}_{\la\mu}=\frac{ 1} {2} \delta^{\a\b}_{\la\mu}:= \frac{ 1} {2} (\delta^\a_\la\delta^\b_\mu-\delta^\a_\mu\delta^\b_\la),\,\,\,\, \forall y\in U_x
  \end{equation}
where $\delta^{\a\b}_{\la\mu}$ is the generalized Kronecker delta. 

\section{Homothetic vector fields that are gradients}\label{sec33}
In this appendix vector fields $X$ that satisfy $\n X=c\mathds{1}$, with $c\in \mathbb{R}$, are studied. They are of particular interest because, on the one hand they are homothetic so that $\pounds_X R=0$ and, on the other hand, they arise naturally is semi-Riemannian manifolds that contain a second-order parallel tensor field $T$ ($\n\n T=0$), see next Appendix \ref{nnT=0}. Observe that  if $c=0$ then $X$ is a parallel vector field. Therefore, if $c=0$ either there is a parallel lightlike vector field or $(M,g)$ is non-degenerately reducible (or decomposable), and the metric decomposes into  $g=\ds\frac{1}{g(X,X)} \uwidehat{X}\otimes \uwidehat{X} + (g-\frac{1}{g(X,X)} \uwidehat{X}\otimes \uwidehat{X} ).$

Consider then the case  when $c\neq 0$.
\begin{lema}\label{hv1}
  Let $(M,g)$ be a semi-Riemannian manifold and $X$ a vector field such that 
  \be\label{nX}
  \n X=c\mathds{1}
  \ee
   for some $c\in\R-\{0\}$. Then:
  \begin{enumerate}[(a)]
    \item $g(X,X)$ is not constant. In particular, $X$ is not lightlike.
    \item $\uwidehat{X}$ is a closed one-form. 
        \item $X$ is homothetic.
        \item $ \n\n T(X,Y)= \n\n T(Y,X),$ for any tensor field $T$ and any vector field $Y$, where $\n\n T(X,Y):=\n_X (\n_Y T)-\n_{\n_X Y}T$.
           \end{enumerate}
   \end{lema}
   
   {\bf Proof.}
   \begin{enumerate}[(a)]
     \item If $g(X,X)=$constant, by taking the covariant derivative here and using \eqref{nX} one would get $2c\uwidehat{X}=0$ and, as $c\neq 0$, this is impossible due to \eqref{nX}.
     \item Usign the covariant version of $\eqref{nX}$ (that is, $\n \uwidehat{X} =c g$), one has $(d\uwidehat{X})_{\a\b}=\n_{[\a}\uwidehat{X}_{\b]}=cg_{[\a\b]}=0.$ 
     \item Since $\n \uwidehat{X}=cg$, if $\pounds_X$ denotes the Lie derivative with respect to $X$ one has that $\pounds_X g (Y,Z) =\n \uwidehat{X} (Y,Z) +\n \uwidehat{X} (Z,Y)=2c g(Y,Z)$ for all vector fields $Y,Z$.
 \item  Since $\n \uwidehat{X}=cg$ one gets that $\n\n \uwidehat{X} =0$ which implies via the Ricci identity
 \be\label{XR}
 X_{\rho} R^{\rho}\,_{\b\la\mu}=0.
 \ee
 Inserting this in the expression of $X^\rho \n_{[\rho}\n_{\b]} T$ via the Ricci identity for any tensor field $T$, one readily gets that $X^\rho \n_{[\rho}\n_{\b]} T=0$ and the assertion follows.
 \end{enumerate}
   \begin{flushright}
{$\blacksquare$}\end{flushright} 
As is well known, $X$ being homothetic it leaves the connection and the curvature invariant
$$\pounds_X \n =0, \hspace{1cm} \pounds_X R=0.$$


\begin{lema}\label{hv3}
    Let $(M,g)$ be a semi-Riemannian manifold and $ X$ a vector field that satisfies \eqref{nX} for some $c\in\R-\{0\}$. Then,
    $$\n_X R+2cR=0,$$
    $$\n_X\n^r R+ c(r+2)\n^r R=0, \hspace{5mm} \forall r\in \mathbb{N}.$$
    \end{lema}
     {\bf Proof.}
     Writing $\mathcal{L}_X R=0$ in abstract index notation:
     $$X^\rho (\n_\rho R^{\a}\,_{\b\la\mu})-(\n_\rho X_\a)R^{\rho}\,_{\b\la\mu}+(\n_\b X^\rho)R^{\a}\,_{\rho\la\mu}+(\n_\la X^\rho)R^{\a}\,_{\b\rho\mu}+(\n_\mu X^\rho)R^{\a}\,_{\b\la\rho}=0$$
  and using here \eqref{nX} one derives
           $$X^\rho (\n_\rho R^{\a}\,_{\b\la\mu})+2cR^{\a}\,_{\b\la\mu}=0.$$          
 Covariantly differentiating this expression and using \eqref{nX} again
         $$X^\rho (\n_\sigma \n_\rho R^{\a}\,_{\b\la\mu})+3c\n_\sigma R^{\a}\,_{\b\la\mu}=0$$     
      and using now the property in Lemma \ref{hv1}, item (d), the previous expression becomes 
      $$X^\rho (\n_\rho \n_\sigma R^{\a}\,_{\b\la\mu})+3c\n_\sigma R^{\a}\,_{\b\la\mu}=0.$$  
         And so on.
           \begin{flushright}
{$\blacksquare$}\end{flushright}

\section{Semi-Riemannian manifolds with a tensor field $T$ such that $\nabla^rT =0$}\label{nnT=0}
In this appendix the question of when there can be tensor fields $T$ with a vanishing $r^{th}$-order covariant derivative is tackled. Apart from its intrinsic interest, this has obvious applications to the case of $r^{th}$-order symmetric semi-Riemannian manifolds, but it has also relevance because the existence of such $T$ implies the existence of homothetic vector fields satisfying \eqref{nX}.

\begin{lema} \cite{TA,Jaen,SN2}\label{nnT}
If an irreducible semi-Riemannian manifold $(M,g)$ carries a non-zero tensor field $T$ that satifies $\n^2 T=0$, then $\ds{X:=\frac{1}{2}{\rm grad} (g(T,T))}$ either 
\begin{enumerate}[(i)]
\item is a vector field that satisfies \eqref{nX} with $c\neq 0$, or
\item is a parallel lightlike vector field and $g(\n T,\n T)=0$, or
\item vanishes and $g(T,T)=$constant.
\end{enumerate}
\end{lema}
{\bf Proof.}
Define the function $f$ as half the total contraction of $T$ with itself
$$ f:=\frac{1}{2} g(T,T)$$
and set $X:= {\rm grad} f$. Due to the assumption $\n^2 T=0$ one has for all vector fields $Y,Z$
$$
{\rm Hess} f (Y,Z) = g(\n_Y T, \n_Z T)
$$
and also 
$$
\n {\rm Hess} f =0.
$$
Hence, Hess$f$ is a 2-covariant symmetric parallel tensor field. Since $(M,g)$ is irreducible a classical theorem due to Eisenhart \cite{Eis} implies that Hess$f = cg$ for some constant $c$ so that \eqref{nX} holds. If $c\neq 0$ then case $(i)$ follows. If $c=0$ then $\n X=0$ and $X$ is parallel and necessarily null due to irreducibility, or zero. The former case is $(ii)$ and the latter is $(iii)$.
 \begin{flushright}
{$\blacksquare$}\end{flushright}

The classical result for positive definite metrics can then be easily deduced from this Lemma.
\begin{theorem}
Let $(M,g)$ be a Riemannian manifold with a tensor field $T$ that satisfies $\n^r T=0$. If either
\begin{enumerate}[(1)]
\item \cite{NO} $(M,g)$ is complete and irreducible
\item \cite{TA} $T$ is any of Riemann, Ricci or Weyl tensors, 
\end{enumerate}
then $\n T=0$.
\end{theorem}
{\bf Proof.}
It is enough to prove the result for $r=2$. To prove $(1)$, notice that either $(i)$ or $(iii)$ in Lemma \ref{nnT} hold. But if an irreducible Riemannian manifold is complete, all homothetic vector fields are in fact Killing vector fields, and thus only $(iii)$ is possible and $\n T=0$.

To prove $(2)$, due to the de Rham decomposition theorem \cite{RH} irreducibility can be assumed and $T$ ---which is either $R$, $Ric$ or $C$--- vanishes in the flat part of the de Rham decomposition. Again from Lemma \ref{nnT} $\n X =c \mathds{1}$ and, if $c\neq 0$, the formulas in Lemma \ref{hv3} hold. In particular $\n_X \n R =0 = -3 c R$ so that $R=0$. Similarly for $Ric$ or $C$. If on the other hand $c=0$, then $(iii)$ in Lemma \ref{nnT} applies and the result also follows. 
 \begin{flushright}
{$\blacksquare$}\end{flushright}

As noticed by S\'anchez \cite{San}, for irreducible Riemannian $(M,g)$ the only possibility left open for having $\n^2 T=0$ without $\n T=0$ is when the manifold is incomplete and has a proper (non-Killing) homothetic vector field of type \eqref{nX}.

One can however prove a fully general theorem if there is at least one generic point in the manifold.
  \begin{theorem}[\cite{SN1}]\label{rsls}
  If a semi-Riemmanian manifold $(M,g)$ has a generic point $x\in M$ and a tensor field $T$ on the manifold that satisfies $\n^r T=0$, then $\n T=0$ everywhere.
\end{theorem}
{\bf Proof.}
It is enough to prove this for $r=2$: $\nabla_{\lambda}\nabla_{\mu}T_{\alpha_{1}\dots \alpha_{q}}=0$.
From the Ricci Identity applied to $\nabla_{[\lambda}\nabla_{\mu]}T_{\alpha_{1}\dots \alpha_{q}}=0$ the  relation
\be\label{step1}
\sum_{i=1}^{q}R^{\rho}{}_{\alpha_{i}\lambda\mu}T_{\alpha_{1}\dots\alpha_{i-1}\rho\alpha_{i+1}\dots \alpha_{q}}=0
\ee
follows. Multiplying by the inverse ${\cal R}^{-1}$ at $x$:
$$
\left. \sum_{i=1}^{q}\left(g_{\alpha_i\lambda}
T_{\alpha_1\dots\alpha_{i-1}\mu\alpha_{i+1}\dots\alpha_q}-
g_{\alpha_i\mu}T_{\alpha_1\dots\alpha_{i-1}\lambda\alpha_{i+1}\dots\alpha_q}\right)\right|_{x}=0
$$
Differentiating \eqref{step1}, then multiplying by ${\cal R}^{-1}$ and using this last expression one also gets at $x$
$$
\left. \sum_{i=1}^{q}\left(g_{\alpha_i\lambda}
\nabla_{\nu}T_{\alpha_1\dots\alpha_{i-1}\mu\alpha_{i+1}\dots\alpha_q}-
g_{\alpha_i\mu}\nabla_{\nu}
T_{\alpha_1\dots\alpha_{i-1}\lambda\alpha_{i+1}\dots\alpha_q}\right)\right|_{x}=0 .
$$
The identity $\nabla_{[\lambda}\nabla_{\mu]}\nabla_{\nu}T_{\alpha_1\dots 
\alpha_q}=0$ provides analogously
\bean
g_{\nu\lambda}\nabla_{\mu}T_{\alpha_1\dots \alpha_q}-
g_{\nu\mu}\nabla_{\lambda}T_{\alpha_1\dots \alpha_q}+
\left. \sum_{i=1}^{q}\left(g_{\alpha_i\lambda}
\nabla_{\nu}T_{\alpha_1\dots\alpha_{i-1}\mu\alpha_{i+1}\dots\alpha_q}-
g_{\alpha_i\mu}\nabla_{\nu}
T_{\alpha_1\dots\alpha_{i-1}\lambda\alpha_{i+1}\dots\alpha_q}\right)\right|_{x}=0.
\eean
Substracting the last two expressions one obtains
$$
\left.g_{\nu\lambda}\nabla_{\mu}T_{\alpha_1\dots \alpha_q}-
g_{\nu\mu}\nabla_{\lambda}T_{\alpha_1\dots \alpha_q}\right|_{x}=0 .
$$
Contracting here $\nu$ and $\lambda$ one derives 
$\nabla_{\mu}T_{\alpha_1\dots \alpha_q}|_{x}=0$.
Finally, as $\nabla T$ is parallel then $\nabla T=0$ on the entire $M$.

\begin{flushright}
{$\blacksquare$}\end{flushright}






%
%
\end{document}